\documentclass{amsart}

\calclayout

\usepackage[unicode]{hyperref}
\usepackage[lite]{amsrefs}
\usepackage{graphicx}
\usepackage{mathtools}
\usepackage{enumitem}
\usepackage{amssymb}

\usepackage{tikz-cd}
\usetikzlibrary{calc}

\newcounter{figuregroup}
\newenvironment{figuregroup}{%
  \global\let\thefiguresave=\thefigure
  \global\let\theHfiguresave=\theHfigure
  \numberwithin{figure}{figuregroup}
  \addtocounter{figure}{1}
  \setcounter{figuregroup}{\value{figure}}
  \setcounter{figure}{0}
}{%
  \global\let\thefigure=\thefiguresave
  \global\let\theHfigure=\theHfiguresave
  \setcounter{figure}{\value{figuregroup}}
}

\newcommand{\uss}[2]{\underset{{\scriptstyle #2}}{#1}}
\newcommand{\supdelta}{\textsuperscript{$\delta$}}
\newcommand{\power}{\mathcal{P}}
\newcommand{\scrF}{\mathcal{F}}
\DeclareMathOperator{\im}{im}
\DeclareMathOperator{\id}{id}
\DeclareMathOperator{\dom}{dom}
\newcommand{\mvert}{\,|\,}      

\newcommand{\trfootnote}[1]{\footnote{#1}}
\newcommand{\bugfootnote}[2]{\footnote{#1[Translator's note #2]}}
\newcommand{\iitem}[1]{\\[0.5\baselineskip]#1\quad}
\newcommand{\lastiitem}{\\[0.5\baselineskip]}

\begin{document}

\title{$S$-glued Sums of Lattices}
\def\translatedby{(translated from the German [T1] by Dale R. Worley;
version of December 10, 2024)}
\author[Ch.\ Herrmann]{Christian Herrmann\\\protect\translatedby}

\maketitle

\begin{quotation}
[\textsc{Translator's Abstract.}
For many equation-theoretical questions about modular lattices, 
Hall and Dilworth give a useful construction:
Let $L_0$ be a lattice with largest element $u_0$, $L_1$ be a lattice
disjoint from $L_0$ with smallest element $v_1$, and $a \in L_0$,
$b \in L_1$
such that the intervals $[a, u_0]$ and $[v_1, b]$ are isomorphic.
Then, after identifying those intervals you obtain $L_0 \cup L_1$, a
lattice structure whose partial order is the transitive relation
generated by the partial orders of $L_0$ and $L_1$. It is modular if $L_0$ and
$L_1$ are modular.
Since in this construction the index set $\{0, 1\}$ is essentially a
chain, this work presents a method --- termed S-glued --- whereby a
general family $L_x\ (x \in S)$
of lattices can specify a lattice with the small-scale lattice
structure determined by the $L_x$ and the large-scale structure
determined by $S$.
A crucial application is representing finite-length modular lattices
using projective geometries.]
\end{quotation}

\textit{The translator welcomes suggestions for improving the translation.
A later exposition of this subject by the same author is [T2].
The author (private communication) recommends Bandelt's exposition [T4] as
an easier approch to the decomposition.
All footnotes are the translator's.
The translator has attempted to preserve the typography of the
original.
Page markings of the form ``\emph{[page xxx/yyy]}'' give first the page of
the PDF of the original, \emph{xxx}, and then the page of the original
journal volume, \emph{yyy}.
Typographical errors have been footnoted and corrected.
In some places, more detail has been provided in proofs.
In three places there are minor mathematical errors; they have been
footnoted and detailed discussion of them is in the Translator's Notes
at the end.  The
errors do not weaken the article's results.}

\hskip 0pt\\{[page 1/255]}\nopagebreak   

\vskip\baselineskip
\centerline{\textbf{0. Introduction}}
\nopagebreak\vskip\baselineskip

For many equation-theoretical questions about modular lattices,
or such questions as arise from embedding problems,
a construction of lattices found to be useful is that of Hall and
Dilworth that was given in [4] (cf. J\'onsson [6, 7] and Gr\"atzer,
J\'onsson and Lakser [3])%
\trfootnote{Day and Freese [T5] suggest that Dilworth's earlier
article [T6] should be
considered the first exposition of lattice gluing.
Dilworth's publications all concern the gluing of two lattices, where
the skeleton (in the terminology of this article) is the two-element
chain, $\textbf{2}$.
Personal accounts from Dilworth's students suggest he was aware of the
possibility of generalizing the skeleton, but that he did not develop it.}%
: Let $L_0$ be a lattice with largest element $u_0$,
$L_1$ be a lattice disjoint from $L_0$ with smallest element $v_1$, and
$a \in L_0$,
$b \in L_1$ elements such that the intervals $[a, u_0]$ and $[v_1, b]$ are
isomorphic to each other. Then, after identifying those intervals
you obtain $L_0 \cup L_1$, a lattice
structure whose partial order is the transitive relation generated by
the partial orders of
$L_0$ and $L_1$. It is modular if $L_0$ and $L_1$ are modular.

Since in this construction the index set $\{0, 1\}$ is essentially
a chain, it makes sense to ask whether a more general
family $L_x\ (x \in S)$ of lattices can specify a lattice
with the small-scale lattice structure determined by the $L_x$
and the large-scale structure determined by $S$.
One approach to this, which particularly applies to modular lattices
and generalizes the Hall--Dilworth construction,
is termed \textit{$S$-glued} in this work.
A significantly different approach
can be found in the $L$-sum of Koh [9], which when applied to
modular lattices does not generally produce modular lattices.

A crucial application is
representing finite-length modular lattices using projective geometries,
which is based on the main theorem to be proven in section~6
and Birkhoff's theorem (cf.\ [1; Chap.~IV, \S~7]) that the atomistic
finite-length modular lattices are precisely the
finite-dimensional (possibly reducible) projective geometries.

\textbf{Main Theorem.} \textit{Every finite-length modular lattice $M$
is an $S$-glued
sum of its maximal atomistic intervals. The lattice $S$ can be
chosen as the set $S(M)$ of the smallest elements of these intervals,
which is a sub-join-semilattice of $M$.}

\hskip 0pt\\{[page 2/256]}\\\nopagebreak

In particular, arithmetic in $M$ is determined by the statements
(10)--(16) --- and the dual statements --- stated in section~2. The lattice
$S(M)$ is called the \emph{skeleton} of $M$.

\vskip\baselineskip
\centerline{\textbf{1. $S$-glued systems of lattices}\trfootnote{We
translate ``verklebe'' as ``glued'' and ``verbundene'' as
``connected''.}}
\nopagebreak\vskip\baselineskip

First of all we want to define that a partially-ordered
set $P$ is \emph{of finite length} if every chain in $P$ is finite, and
\emph{of finite length $l(P)$}, if
$l(P) = \sup \{ |K|- 1 \mvert K \textrm{ is a chain in } P\} < \infty$.
Furthermore, $a \prec b$ means that $a$ is a lower neighbor of $b$.

Now let $S = (S, \vee, \wedge)$ be a finite-length lattice and
$L_x = (L_x, \uss +x, \uss \cdot x)\ (x \in S)$
be a family of finite-length lattices.  We for now
assume that the underlying sets of the $L_x$ have already been chosen such that
the construction promised in the introduction can be made without
having to either further identify elements or make them disjoint
--- the prerequisites for this are incorporated into the concept of
$S$-glued system of lattices. Later (in section~4)
will be shown how to produce such $S$-glued systems.

A family $L_x\ (x \in S)$ of lattices indexed by a lattice $S$
is called an \emph{$S$-glued
system}, provided that for all $x, y \in S$:
\iitem{(1)}\textit{If $x \leqq y$ and $L_x \cap L_y \neq \emptyset$, then
$L_x \cap L_y$ is a filter of $L_x$ and
an ideal of $L_y$.}
\iitem{(2)}\textit{Under the assumptions of (1), for all $a, b \in L_x \cap L_y$,
$a \leqq b$ in $L_x$ if and only if $a \leqq b$ in $L_y$.}
\iitem{(3)}\textit{If $x \prec y$, then $L_x \cap L_y \neq \emptyset$.}
\iitem{(4)}\textit{$L_x \cap L_y \subseteq L_{x \wedge y} \cap L_{x \vee y}$.}
\lastiitem
The purpose of conditions (1) and (3) is to identify
duplicate elements,
while (1) and (4) together prevent undesired collision of elements.
The following families $L_x\ (x \in S)$ are not $S$-glued systems ---
where differently labeled elements are assumed to be different:
\\
\includegraphics[width=\linewidth]{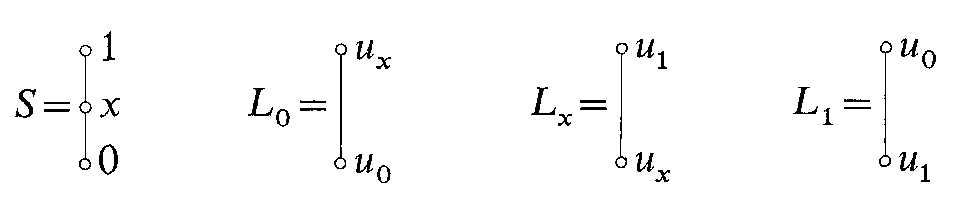}
\\
and also
\\
\includegraphics[width=\linewidth]{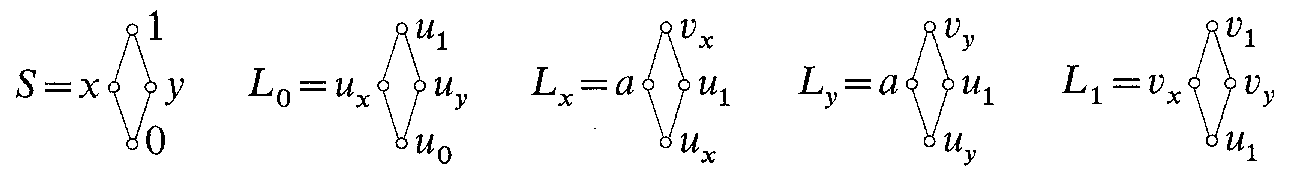}
\\[\baselineskip]{[page 3/257]}\\[\baselineskip]\nopagebreak
since in the first example $L_0 \cap L_1 = \{u_0\}$ is not a filter of
$L_0$ (i.e.\ (1)
is violated) and in the second $L_x \cap L_y = \{u_1,
a\} \nsubseteq L_0 \cap L_1$, so (4) is violated, even though
(1)--(3) are satisfied.

The $\uss\leqq x$ are the partial orders of the lattices $L_x$, $0_x$
are the smallest and $1_x$ the largest elements of $L_x$. For
$a, b \in L_x$ let $[a,b]_x$ be the interval in $L_x$
with smallest element $a$ and largest element $b$. With these definitions
we get the following simple conclusions from the axioms (1)--(4):
\iitem{(5)}\textit{If $x \leqq y$ and $L_x \cap L_y \neq \emptyset$, then $0_y,
1_x \in L_x \cap L_y$ and $L_x \cap L_y = [0_y, 1_x]_x = [0_y, 1_x]_y$.}
\iitem{(6)}\textit{If $x \leqq z \leqq y$, then $L_x \cap L_y \subseteq L_z$.}
\iitem{(7)}\textit{$L_x \cap L_y = L_{x \wedge y} \cap L_{x \vee y}$.}
\iitem{(8)}\textit{\emph{[added by tr.:} If $L_x \cap L_y \neq \emptyset$,\emph{]}
$L_x \cap L_y$ is an interval in $L_x$ and in $L_y$ and
for all $a,b \in L_x \cap L_y$,
$a \uss\leqq x b$ if and only if $a \uss\leqq y b$; also $a \uss+x
b=a \uss+y b$ and $a \uss\cdot x b=a \uss\cdot y b$.}
\iitem{(9)}\textit{If $a \uss\leqq x b$ and $a \in L_y$, then
$a \uss\leqq{x \vee y} b$.}
\lastiitem

\textit{Proof of (5)--(9).} Regarding (5): Let $a \in L_x \cap
L_y$. Then $0_y \uss\leqq y a$, so
by (1), $0_y \in L_x$ and by (2), $0_y \uss\leqq x
a$. Dually $a \uss\leqq x 1_x$, $1_x \in L_y$ and $a \uss\leqq y 1_x$.
Hence $L_x \cap L_y \subseteq [0_y, 1_x]_x \cap [0_y,
1_x]_y$. Conversely, if $0_y \uss\leqq x a \uss\leqq x 1_x$, then
$a \in L_x \cap L_y$ due to (1). Thus $L_x \cap L_y = [0_y, 1_x]_x$,
Likewise, we get $L_x \cap L_y = [0_y, 1_x]_y$.

Regarding (6): First, let $x \prec z \leqq y$. If $L_x \cap L_y
= \emptyset$ then nothing need be proven. Otherwise by (5), $L_x
\cap L_y = [0_y, 1_x]_x$ and by (3),
$L_x \cap L_z = [0_z, 1_x]_x$. This means that $L_z \cap L_y$ is not empty,
which --- again because of (5) -- means $0_y \in L_z$. So
$0_z \uss\leqq z 0_y$ and
$0_z \uss\leqq x 0_y$, which is why $L_x \cap L_y = [0_y,
1_x]_x \subseteq [0_z, 1_x]_x \subseteq L_z$. The general case
now follows by induction on $z$ in $S$:
Let $x \leqq z^\prime \prec z \leqq y$; by the induction
hypothesis $L_x \cap L_y \subseteq\mskip -\thickmuskip%
\trfootnote{The original has the typo $\leqq$ for $\subseteq$.}\;
L_{z^\prime}$, and, as just shown,
$L_{z^\prime} \cap L_y \subseteq L_z$; together showing
$L_x \cap L_y \subseteq L_{z^\prime} \cap L_y \subseteq L_z$.
(7) is a direct consequence of (4) and (6).

Regarding (8):
\emph{[added by tr.:}  By (7), $L_x \cap L_y = L_{x \vee y} \cap L_{x \wedge y}$.
By (1) and (2), the smallest element of this intersection is $0_{x \vee y}$
and the largest element is $1_{x \wedge y}$, under all of the orders
$\uss\leqq x$, $\uss \leqq y$, $\uss \leqq {x \vee y}$, and $\uss \leqq
{x \wedge y}$.\emph{]}
Let $a, b \in L_x \cap L_y$. Because of (4), $a, b \in
L_x \cap L_y \cap L_{x \vee y}$. By (2),
$a \uss\leqq x b$ is equivalent to $a \uss\leqq{x \vee y} b$, which in
turn is equivalent to $a \uss\leqq y b$, and vice-versa.
Therefore $a \uss\leqq x b$ if and only if $a \uss\leqq y b$.
Let $c \in L_x$ with $a \uss\leqq x c \uss\leqq x b$.
Then by (1), $c \in L_x \cap L_{x \vee y}$ and dually $c \in
L_x \cap L_{x \wedge y}$, so $c \in L_x \cap L_y$ by
(7).
\emph{[added by tr.:}
This shows that $L_x \cap L_y$ is convex in $L_x$, and having both a
smallest and greatest element,\emph{]}
it is an interval in $L_x$ and similarly also in $L_y$.
From this it finally follows that
\begin{equation*}
a \uss+x b = \sup_{L_x} (a,b) = \sup_{L_x \cap L_y} (a,b) =
\sup_{L_y} (a, b) = a \uss+y b \textrm{\quad and dually\quad}
a \uss\cdot x b = a \uss\cdot y b.
\end{equation*}

\hskip 0pt\\{[page 4/258]}\\\nopagebreak

Regarding (9): Let $a \uss\leqq x b$ and $a \in L_y$. By (4), we have
$a \in L_{x \vee y}$, so by (1) and (2), $b \in L_{x \vee y}$ and
$a  \uss\leqq{x \vee y} b$.

\vskip\baselineskip
\centerline{\textbf{2. $S$-glued sums of lattices}}
\nopagebreak\vskip\baselineskip

On the (by no means disjoint) union $L := \bigcup_{x\in S} L_x$ of the lattices
of an $S$-glued system, let $\leqq$ be the transitive closure of
the union of the (reflexive and transitive) relations $\uss\leqq x$.
$L = (L, \leqq)$ is
called the \textit{$S$-glued sum} of $L_x\ (x \in S)$.

\textbf{Theorem 2.1.} \textit{The $S$-glued sum $L$ of lattices
$L_x\ (x \in S)$ is a
finite-length lattice. $a \prec b$ holds in $L$ if and only if
$a \uss\prec x b$ for some $x \in S$.
For every $x \in S$ the lattice $L_x = (L_x, \uss+x, \uss\cdot x)$ is
an interval sublattice of $L$.
We write: $a + b := \sup_L (a, b)$ and $a \cdot b := \inf_L (a, b)$.}

\textbf{Addendum\trfootnote{We translate ``Zusatz'' as ``addendum''.} 2.2.}
 \textit{The map $x \mapsto 0_x$ is a
$\sup$-homomorphism of $S$
into $L$, and $x \mapsto 1_x$ is an $\inf$-homomorphism of $S$ into $L$.
If the given $S$-glued lattice system is \emph{monotone}%
\bugfootnote{To
make addendum~2.2 correct, the definition of ``monotone'' must be stricter.
Specifically, an $S$-glued lattice system is \emph{monotone} if for all
$x \prec y$, $x, y \in S$, then $L_y \nsubseteq L_x$ and
$L_x \nsubseteq L_y$.  Note this renders the definition of ``monotone''
self-dual.}{2}
--- $L_y \nsubseteq L_x$ for all $x \prec y$ --- then these maps are
injections.}

\textit{Proof.} The partial-order structure of $L$ is characterized more
precisely (and proven to be a lattice\trfootnote{The original says
``partial order'' but ``lattice'' is clearly needed, and proven.}) by
these statements, proven below:
\iitem{(10)}\textit{$a \leqq b$ holds in $L$ if and only if
there exists $x_1 \leqq x_2 \leqq \cdots \leqq x_n$ in $S$ and
$a_i \in L_{x_i} \cap L_{x_{i+1}} (i=1 \ldots, n-1)$ such that
$a = a_0 \uss\leqq{x_1} a_1 \uss\leqq{x_2} \cdots \uss\leqq{x_n} a_n = b$.
If $a \in L_x$, then we can choose $x_1 \geqq x$.}
\iitem{(11)}\textit{If $a \leqq c \leqq b$ in $L$ and $a, b \in L_x$ for some $x \in S$,
then $c \in L_x$ and $a \uss\leqq x c \uss\leqq x b$.}
\lastiitem

The antisymmetry of the relation $\leqq$ follows directly from (11),
and the finite length of $L$ follows directly from
(10).\bugfootnote{Finite length is not guaranteed without further
conditions.}{1}
Furthermore, (11) says that every
$L_x$ is an interval of $L$ and $a \prec b$ follows from $a \uss\prec x b$
(the converse
of this is trivial). It should be noted that if $S$ has finite length $l$ and
every $L_x$ has finite length $l_x \leqq k$ for a fixed $k$, $L$ has
length at most $k(l + 1)$.

Calculating a supremum in $L$ (and proving its existence)
is done with the help of these statements, proven below:
\iitem{(12)}\textit{If $x \leqq y$, then $0_x \leqq 0_y$;
\emph{[added by tr.:} given $x \leqq y$,\emph{]} $0_x < 0_y$ if and only if
$x < y$ and $L_y \nsubseteq L_x$\bugfootnote{The original is
incorrect here; the condition is $L_x \nsubseteq L_y$.  Note that the
dual statement, (12\textsuperscript{*}) is ``\textit{If $x \leqq y$,
then $1_x \leqq 1_y$;
given $x \leqq y$, $1_x < 1_y$ if and only if
$x < y$ and $L_y \nsubseteq L_x$.}''}{2}.}
\iitem{(13)}$\sup_L (0_x, 0_y) = 0_{x \vee y}$.
\iitem{(14)}\textit{If there is an $x \in S$ with $a, b \in L_x$, then
$\sup_L (a, b)$ exists and $= a \uss+x b$.}
\\[\baselineskip]{[page 5/259]}\\\nopagebreak   
\iitem{(15)}\textit{If $x = x_0 \prec x_1 \prec \cdots \prec x_n = y$ in $S$
with $n \geqq 0$\trfootnote{The original has ``$n \geqq 1$'' but the proof
shows that ``$n \geqq 0$'' was intended.}, $a \in L_x$ and
$b \in L_y$, then}
\begin{equation*}
\sup_L (a, b) \textrm{ exists\quad and\quad}
= ((\cdots((a \uss+{x_0} 0_{x_1}) \uss+{x_1}
0_{x_2})\uss+{x_2} \cdots) \uss+{x_{n-1}} 0_{x_n}) \uss+y b.
\end{equation*}
\iitem{(16)}\textit{If $a \in L_x$ and $b \in L_y$, then}
\begin{equation*}
\sup_L (a, b) \textrm{ exists\quad and\quad}
= \sup_L(a, 0_{x \vee y}) \uss+{x \vee y} \sup_L(b, 0_{x \vee y}).
\end{equation*}
\lastiitem

Infima in $L$ can be determined using the dual statements
(13\textsuperscript{*})-(16\textsuperscript{*}).
Because of (14) and (14\textsuperscript{*}), the $L_x$ are interval sublattices
of $L$. By (13) (resp.\ (13\textsuperscript{*})), the map $x \mapsto 0_x$
(resp.\ $x \mapsto 1_x$) is a $\sup$(resp.\ $\inf$)-homomorphism, and
because of (12),
they are injective if the $S$-glued system is monotone. Thus, to
prove the theorem and corollary it remains only
to show (10)--(16).

\textit{Proof of statements (10)-(16).}  Regarding (10): Let $a \leqq
b$ in $L$. Then there is,
by the definition of $\leqq$, $x_1, \ldots, x_n \in S$ and
$a_i \in L_{x_i} \cap L_{x_{i+1}} (i= 1, \ldots n-1)$
such that $a = a_0 \uss\leqq{x_1} a_1 \uss\leqq{x_2} \cdots \uss\leqq{x_n}
a_n = b$.
Repeated application of (9)
shows that we can choose the $x_i$ such that
$x_i \leqq x_{i+1}\ (i = 1, \ldots, n-1)$ and, if $a \in L_x$, $x \leqq x_1$.
The converse is trivial.

Regarding (11): In view of (10), it is sufficient to show that
$a_0 \uss\leqq{x_1} a_1 \uss\leqq{x_2} \cdots \uss\leqq{x_n} a_n$,
$a_0, a_n \in L_x$ and $x \leqq x_1 \leqq \cdots \leqq x_n$\trfootnote{The
original includes an $x_0$ in this sequence, even though it is never
used.}
implies $a_0 \uss\leqq x a_i \uss\leqq x a_n$
for $i = 0, \ldots, n$. This is done by induction on $n$. The case $n = 0$
is trivial. Now assume $n \geqq 1$.
The induction hypothesis is $a_0 \uss\leqq x a_i \uss\leqq x a_{n-1}$
for all $i = 0, \ldots, n-1$.
Thus $a_{n-1} \in L_x$ and so $a_{n-1} \in L_x \cap L_{x_n}$.
Since $a_{n-1} \uss\leqq{x_n} a_n$ and $x \leqq x_n$,
by (1), $a_{n-1} \uss\leqq{x} a_n$, 
and thus also $a_0 \uss\leqq x a_i \uss\leqq x
a_n$ for $i = 0, \ldots, n$.

Regarding (12): Assume $x \prec y$. Then from (3) and (5) it follows that
$0_y \in L_x$, so that $0_x \uss\leqq x 0_y$ and $0_x \leqq 0_y$.
By induction, for any $x \leqq y$, $0_x \leqq 0_y$.

Again, assume $x \prec y$, so by (3), $L_x \cap L_y \neq \emptyset$.
If $L_y \nsubseteq L_x$, then $0_x \neq 0_y$, since otherwise by (5),
$L_x = [0_x, 1_x]_x = [0_y, 1_x]_x = L_x \cap L_y$, implying
$L_x \subseteq L_y$.
In general, if $x \leqq y$, 
if $L_y \nsubseteq L_x$,
there are $x^\prime, y^\prime$ with $x \leqq x^\prime \prec y^\prime \leqq y$
and $L_{y^\prime} \nsubseteq L_{x^\prime}$,\bugfootnote{To be correct, the
roles of $x$ and $y$ must be swapped in the preceding clause of this
sentence.}{2} and thus $0_x \leqq
0_{x^\prime} < 0_{y^\prime} \leqq 0_y$.

Regarding (13): $0_{x \vee y}$ is an upper bound of $0_x$ and $0_y$
because of (12).
Choose $e \geqq 0_x, 0_y$ in $L$, i.e.
\begin{equation*}
0_x = a_0 \uss\leqq{x_1} a_1 \uss\leqq{x_2} \cdots \uss\leqq{x_n} a_n = e
\textrm{\quad and\quad}
0_y = b_0 \uss\leqq{y_1} b_1 \uss\leqq{y_2} \cdots \uss\leqq{y_m} b_m = e
\end{equation*}
with $x \leqq x_1 \leqq \cdots \leqq x_n$ and $y \leqq y_1 \leqq \cdots \leqq y_m$.
By (4), $e \in L_{x_n \vee y_m}$, therefore
$e \geqq 0_{x_n \vee y_m} \geqq\mskip -\thickmuskip\trfootnote{The
original incorrectly has
$=$ for $\geqq$; $e$, $x_n$, and $y_m$ might be very large.}\;
0_{x \vee y}$ by (12).

Regarding (14): The proof is by descending induction
on $x$ in $S$. It must be shown that for any $e > a, b$, we can
conclude $e \geqq a \uss+x b$.

Let $x = 1$ ($1$ being the largest element of $S$), so $a, b \in L_1$
and $e > a,b$.
$1_1$ is the
largest element of $L$, so $a, b \leqq e \leqq  1_1$ and from this and
(11), $e \in L_1$ and
$a, b \uss\leqq1 e$. Therefore $a \uss+1 b \uss\leqq1 e$ and
$a \uss+1 b \leqq e$.

\hskip 0pt\\{[page 6/260]}\\\nopagebreak

For general $x \in S$, we have the induction hypothesis:
\iitem{(*)}For all $z > x$ and $g, h \in L_z$ we have
$g \uss+z h = \sup_L(g, h)$%
\trfootnote{The original is missing a line break here.}
\lastiitem
and the induction goal:
For all $a, b \in L_x$, $a \uss+x b = \sup_L (a, b)$.

The latter can in turn be proven inductively
by $(a, b)$ descending in $L_x \times L_x$.
The induction start, $a = b = 1_x$, is trivial. So let $a, b \in L_x$
with the induction hypothesis:
\iitem{(**)}For all $c, d \in L_x$ with $c \geqq a$, $d \geqq b$ and
either $c \neq a$ or $d \neq b$,
\begin{equation*}
c \uss+x d = \sup_L(c, d).
\end{equation*}
\lastiitem
Let $e \in L$ and $e > a, b$. Because of (10) this means
\begin{equation*}
a = a_0 \uss<{x_1} a_1 \uss<{x_2} a_2 \uss<{x_3} \cdots \uss<{x_n} a_n = e
\textrm{\quad and\quad}
b = b_0 \uss<{y_1} b_1 \uss<{y_2} b_2 \uss<{y_3} \cdots \uss<{y_m} b_m = e
\end{equation*}
(collapsing successive values that are equal) with $x \leqq x_1 \leqq
x_2 \leqq \cdots \leqq x_n$ and
$x\trfootnote{The original has the typo $y$ for $x$.}
\leqq y_1 \leqq y_2 \leqq \cdots \leqq y_m$.
If $x_1 = x$, then from (**) $a_1 \uss+x b \leqq e$;
but $a \uss+x b \uss\leqq x a_1 \uss+x b$, i.e.\ $a \uss+x b \leqq e$.
Analogously if $y_1 = x$.

The remaining case is $x_1, y_1 > x$.
Since $x \leqq x_1 \wedge y_1 \leqq x_1$ and $a \in L_x, L_{x_1}$,
from (6) it follows that $a \in L_{x_1 \wedge y_1}$.
Similarly, $b \in L_{x_1 \wedge y_1}$.
In the case $x_1 \wedge y_1 > x$ we can therefore use
(*) to conclude that $a \uss+{x_1 \wedge y_1} b \leqq e$.
Then because of (1), $a \uss+x b \in L_x \cap L_{x_1 \wedge y_1}$ and
because of (8), $a \uss+x b = a \uss+{x_1 \wedge y_1} b \leqq e$.

The remaining case is reduced to $x_1, y_1 > x$ and $x_1 \wedge y_1 = x$.
(5) ensures that $0_{x_1}, 0_{y_1} \in L_x$.
Let $c := 0_{x_1} \uss+x 0_{y_1}$. By (1), $c \in L_{x_1}$ and
$c \in L_{y_1}$,
so by (7), $c \in L_{x_1} \cap L_{y_1} =
L_x \cap L_{x_1 \vee y_1}$. Hence $c \geqq 0_{x_1 \vee y_1}$.  On the other hand
by (5), $0_{x_1}, 0_{y_1} \uss\leqq x 0_{x_1 \vee y_1}$ and with (13),
$c = 0_{x_1 \vee y_1} = \sup_L(0_{x_1}, 0_{y_1})$.
By (8) and (*), this results in:
\begin{equation*}
a \uss+x c = a \uss+{x_1} c = \sup_L (a, c) = \sup_L (a, 0_{x_1} , 0_{y_1})
\textrm{\quad and\quad}
b \uss+x c = \sup_L (b, 0_{x_1}, 0_{y_1}).
\end{equation*}
However $0_{x_1} \leqq a_1 \leqq e$, $0_{y_1} \leqq b_1 \leqq e$ and
therefore $a \uss+x c \leqq e$ and $b \uss+x c \leqq e$.

If $a \ngeqq c$, then also $a \uss\ngeqq x c$ and
therefore $a \uss+x c \uss>x a$. (**) guarantees
$(a \uss+x c) \uss+x (b \uss+x c) \leqq e$,
from which immediately follows $a \uss+x b \leqq e$.
We similarly handle when $b \ngeqq c$.

So ultimately the case $a \geqq c$ and $b \geqq c$ remains.
By (12\textsuperscript{*}),
$c \leqq a \leqq 1_x \leqq 1_{x_1} \leqq 1_{x_1 \vee y_1}$ and similarly
for $b$, so
(11) shows $a$ and $b$ are elements of $L_{x_1 \vee y_1}$.
Since $x_1 \vee y_1 > x$, we can
conclude with (*) that $a \uss+{x_1 \vee y_1} b \leqq e$.
However, $a \uss+{x_1 \vee y_1} b = a \uss+x b$ by (8).

Regarding (15): The proof proceeds inductively on $n$.
The case $n = 0$ is direct from (14).
So let $n \geqq 1$.
Set $c_0 := a \in L_{x_0}$.
For $1 \leqq i \leqq n$, we inductively define
$c_i := c_{i-1} \uss+{x_{i-1}} 0_{x_i}$ and prove $c_i \in L_{x_i}$
using (3), (5), and (1), since $0_{x_{i+1}} \in L_{x_i}$.
By the induction on $n$, $c_{n-1} = \sup_L(a, 0_{x_{n-1}})$.
Finally set $d := c_n = c_{n-1} \uss+{x_{n-1}} 0_{x_n}$.
Because of (14)
$d = \sup_L(a, 0_{x_{n-1}}, 0_{x_n}) = \sup_L(a, 0_{x_n})$.
From (2), $d \in L_{x_n}$, so $d \uss+{x_n} b$ is defined
\\[\baselineskip]{[page 7/261]}\\[\baselineskip]\nopagebreak
and by (14), $d \uss+{x_n} b = \sup_L(d, b)$.
Therefore $d \uss+{x_n} b = \sup_L(a,0_{x_n}, b) = \sup_L(a, b)$.

Regarding (16): By (15), $\sup_L(a, 0_{x \vee y}), \sup_L(b,
0_{x \vee y}) \in L_{x \vee y}$; further
\begin{equation*}
\sup_L(a, 0_{x \vee y}) \uss+{x \vee y} \sup_L(b, 0_{x \vee y}) =
\sup_L(a, b, 0_{x \vee y}) = \sup_L(a, b, 0_x, 0_y)= \sup_L(a, b)
\end{equation*}
by (14) and (13).

\vskip\baselineskip
\centerline{\textbf{3. Transfer of lattice-theoretical properties}}
\nopagebreak\vskip\baselineskip

In the this section the concept of the $S$-glued sum is justified
as an aid for the theory of modular lattices
by showing that modularity (but also breadth $\leqq n$,
$n$-distributivity) is transferred from the components to the sum.
Also

\textbf{Lemma 3.1.} \textit{Let $L$ be the $S$-glued sum of $L_x\ (x \in S$)
and $U$ be a subset
of $L$ with smallest element $u$ and largest element $v$.
Let $A := \{ a \mvert a \in U \textrm{ and there is no } b \in U
\textrm{ with } u < b < a \}$ be the set of atoms of $U$.
Finally assume that for all $a \in A$, $u \prec a$ holds in $L$ and that
$v = \sup_L A$.
Then there is an $x \in S$ such that $U \subseteq L_x$ and the lattice structure
of $U$ in $L$ agrees with that in $L_x$, i.e.\ that for all $a, b, c \in U$,
$a + b = c$
($a \cdot b = c$) if and only if $a \uss+x b = c$ $(a \uss\cdot x b = c$).}

\textit{Proof.} By theorem~2.1, for every $a \in A$
there is an $x_a \in S$ such that $u \uss\prec{x_a} a$.
Set $x := \sup_S \{x_a \mvert a \in A\}$.
By (4), $u \in L_x$ and by (1), $a \in L_x$ for
all $a \in A$. Thus $v = \sup_{L_x} A \in L_x$ and
$U \subseteq [u, v] = [u, v]_x \subseteq L_x$.
Now, the agreement of the lattice structures stated above results
directly from theorem~2.1.

\textbf{Theorem 3.2.} \textit{The $S$-glued sum of modular
(semi-modular) lattices is modular (semi-modular).}

\textit{Proof.} Let $L$ be the $S$-glued sum of semi-modular lattices
$L_x\ (x \in S$). If $a \prec b, c$ in $L$ and $b \neq c$,
then $a \uss\prec x b, c$
holds for some $x \in S$ by the lemma (with
$U = \{ a, b, c, b + c\}$) and thus because of the semi-modularity of $L_x$,
$b + c = b \uss+x c \uss\succ x b, c$;
so by theorem~2.1, $b + c \succ b, c$.
Therefore $L$ is of finite length and semi-modular.
In addition, $L$ satisfies
the Jordan-Dedekind chain condition (the proof in [1; p.~40,
Thm.~14] can be generalized suitably).
Similarly, if the $L_x$ are dual semi-modular, $L$ is dual
semi-modular.
In the case of modularity we use the fact that
that for finite lattices, semi-modularity and dual semi-modularity
modularity taken together are equivalent to modularity
(cf. [1; p.~41, Thm.~16]).

\hskip 0pt\\{[page 8/262]}\\\nopagebreak

It is obvious that length or width, for example, are not usually
transferred to the $S$-glued sum.  For breadth and $n$-distributivity,
however, we can prove transfer theorems that enable our later
discussion regarding the ``inner dimension'' of finite-length modular
lattices.  Recall these definitions:

The \textit{breadth} $b(L)$ of a lattice $L$ is the supremum of all
$n$ for which there is a mapping $\phi$ of the Boolean lattice
$\boldsymbol{2}^n$ into $L$ such that for all $a,
b \in \boldsymbol{2}^n$, $\phi\,a \leqq \phi\,b$ if and only if
$a \leqq b$ (this is, $\phi$ is an isomorphism of partial orders).
By A. Huhn [5], a lattice is called \textit{$n$-distributive} if it
is modular and the equation
\begin{equation*}
x \cdot \sum_{i=0}^n y_i =
\sum_{j=0}^n \left( x \cdot \sum_{i=0, i\neq j}^n
y_i \right)\trfootnote{The original has the typo
$\prod_{j=0}^n \left( x + \prod_{i=0, i\neq j}^n y_i \right)$
for the right-hand side.}
\end{equation*}
holds, or the equivalent condition that there is no
sublattice $U$ of $L$ isomorphic with
$\boldsymbol{2}^{n+1}$
with atoms $a_0, \ldots, a_n$ and
an element $w$ of $L$ such that $a_i \cdot w = \inf_L U$
and $a_i + w = \sup_L U$ for $i = 0, \ldots, n$.

\textbf{Theorem 3.3.} \textit{The $S$-glued sum of semi-modular
lattices of breadth $\leqq n$ is semi-modular with breadth $\leqq n$.}

\textit{Proof.} Let $L$ be the $S$-glued sum (semi-modular according
to theorem~3.2) of
semi-modular lattices $L_x\ (x \in S)$ of breadth $\leqq n$. Now let be $U$ a
partially-ordered subset of $L$ isomorphic to $\boldsymbol{2}^m$
with coatoms $b_0, \ldots, b_{m-1}$,
smallest element $u$ and largest element $v$.
Without loss of generality,
we can assume that $U$ is a sub-$\cdot$-semilattice of $L$.
Induction over the length $l$ of the interval $[u, v]$
shows that there is an $x \in S$ with $U \subseteq L_x$.
This is trivial in the case $l=0$.
Assume $l > 0$.
If the atoms $a_i = \prod^{m-1}_{j=0, j\neq i} b_j$ of $U$ are all
upper neighbors of $u$ in $L$,
then the sub-$+$-semilattice of $U$ generated by them (and also isomorphic
to $\boldsymbol{2}^m$) is, by lemma~3.1 a sub-$+$-semilattice
of $L$ and contained entirely in a single $L_x$.
However, if there is an element $a \in L$ with, for example, $u \prec a < a_0$,
so $a \nleqq b_0$ and because of semi-modularity $a_0 \nleqq a + b_0$.
So $a + b_0, b_1, \ldots, b_{m-1}$ generates a sub-$\cdot$-semilattice
of L that is isomorphic to $\boldsymbol{2}^m$,
that contained in the interval $[a, v]$, has length $< l$, and thus by
the induction hypothesis is contained in some $L_x$.
This proves $m \leqq n$ is proved, and so the sum has breadth $\leqq n$.

\textbf{Theorem 3.4.} \textit{The $S$-glued sum of $n$-distributive
lattices is $n$-distributive.}

\textit{Proof.} Let $L$ be the $S$-glued sum (modular by theorem~3.2)
of $n$-distributive lattices $L_x\ (x \in S)$.  Suppose $a_0, \ldots,
a_{m-1}$ are the atoms of a sublattice of $L$ isomorphic to $\boldsymbol{2}^m$,
with smallest element $u$ and largest element $v$, and finally let $w$
be an element of $L$ with $w \cdot a_i = u$ and
\\[\baselineskip]{[page 9/263]}\\[\baselineskip]\nopagebreak
$w + a_i = v$ for $i=0, \ldots, m-1$. Now let $w^\prime \in L$
with $w \leqq w^\prime \prec v$. Define
inductively: $c_0 := w^\prime \cdot a_0$, $a_i^0 := a_i + c_0$
for $i = 0, \ldots, m - 1$; and for $k = 1, \ldots, m-1$,
$c_k := w^\prime \cdot (a_k +c_{k-1})$ and
$a^k_i := a_i + c_k$ for $i = 0, \ldots, m - 1$. Then follows
inductively from modularity:
\begin{equation*}
c_k = \sum_{j=0}^k w^\prime \cdot a_j \leqq c_{k+1} < a^k_{k+1} \leqq v
\end{equation*}
and $a_0^{k+1}, \ldots, a_{m-1}^{k+1}$
are atoms of a sublattice of $L$ isomorphic to $\boldsymbol{2}^m$
with smallest
element $c_{k+1}$ and largest element $v$.
Finally, $w^\prime + a_i^{m-1} = v$ and
$w^\prime \cdot a_i^{m-1} = c_{m-1}$ for $i = 0, \ldots, m - 1$.
Therefore the $a_i^{m-1}$ are upper neighbors
of $c_{m-1}$ in $L$ and by lemma~3.1, $[c_{m-1}, v] \subseteq L_x$
for some $x \in S$.
So it follows from the $n$-distributivity of $L_x$ that $m \leqq n$.

This proof is a modification of an idea of Huhn's.

\vskip\baselineskip
\centerline{\textbf{4. $S$-connected systems}}
\nopagebreak\vskip\baselineskip

We now turn to the question posed in section~1, when
from a given family of lattices and a system
of ``identifications'' between these lattices an $S$-glued
system can be constructed.

Again, let $S$ be a finite-length lattice and $L_x\ (x \in S)$ be a family
of finite-length lattices. Furthermore, let
$\phi_{yx}\ (x, y \in S, x \leqq y)$ be a family
of bijective mappings. The system consisting of the $L_x\ (x \in S)$
and the $\phi_{yx}\ (x, y \in S, x \leqq y)$ is called an
\textit{$S$-connected system} if for all $x, y \in S$:
\iitem{(17)}\textit{If $\phi_{yx} \neq \emptyset$, then $\phi_{yx}$ is an isomorphism
from a filter of $L_x$ to
an ideal of $L_y$; $\phi_{xx}$ is the identity map of $L_x$.}
\iitem{(18)}\textit{If $x \prec y$ in $S$, then $\phi_{yx} \neq \emptyset$.}
\iitem{(19)}\textit{For every $z \in S$ with $x \leqq z \leqq y$ we have
$\phi_{yx} = \phi_{yz} \circ \phi_{zx}$.}
\iitem{(20)}$\im \phi_{(x \vee y)x} \cap \im \phi_{(x \vee y)y} \subseteq
\im \phi_{(x \vee y)(x \wedge y)}$.
\iitem{(20\supdelta)}$\dom \phi_{x(x \wedge y)} \cap
\dom \phi_{y(x \wedge y)} \subseteq \dom \phi_{(x \vee y)(x \wedge y)}$.
\lastiitem

For a mapping $\phi$, $\dom \phi$ always denotes its domain
and $\im \phi$ its range.

Now let $L$ be the disjoint union of the $L_x\ (x \in S)$ and for each
$x \in S$, $i_x$ is the canonical embedding of $L_x$ in $L$.
Then the lattices $K_x := i_x(L_x)$ and the maps $\psi_{yx} := i_y \circ
\phi_{yx} \circ i_x^{-1}$ trivially form an $S$-connected system.

\textbf{Proposition\trfootnote{We translate ``Hilfssatz'' as
``proposition''.} 4.1.} \textit{Let $a \in K_x$ and $b \in K_y$.
Then the following are equivalent:}
\begin{enumerate}
\item[(i)] \textit{There is a $z \in S$ with $\psi_{zx}\,a =
\psi_{zy}\,b$\trfootnote{The original has the typo $\psi_{yz}\,b$ for
for $\psi_{zy}\,b$.}.}
\end{enumerate}
\vskip\baselineskip{[page 10/264]}\\\nopagebreak   
\begin{enumerate}
\item[(ii)] \textit{$\psi_{(x \vee y)x}\,a =
\psi_{(x \vee y)y}\,b$\trfootnote{The original has the typo $\psi_{(x \vee y)x}\,b$
for $\psi_{(x \vee y)y}\,b$.}.}
\item[(iii)] \textit{There is a $z \in S$ with $\psi_{xz}^{-1}\,a
= \psi_{yz}^{-1}\,b$.}
\item[(iv)\footnotemark]
\textit{$\psi_{x(x \wedge y)}^{-1}\,a = \psi_{y(x \wedge y)}^{-1}\,b$.}
\end{enumerate}%
\footnotetext{The original has the typo ``(vi)'' for ``(iv)''.}

\textit{Proof.} The equivalence of (i) and (ii) or (iii) and (iv) is
immediately seen from (19).  Let $a \in K_x$, $b \in K_y$ and
let $\psi_{(x \vee y)x}\,a = \psi_{(x \vee y)y}\,b$.
Then because of (20)\trfootnote{The original has the typo
(20\supdelta) for (20).} there is a
$c \in K_{x \wedge y}$ with $\psi_{(x \vee y)(x \wedge y)}\,c =
\psi_{(x \vee y)x}\,a$.
Because of (19),
$\psi_{(x \vee y)(x \wedge y)}\,c = \psi_{(x \vee
y)x}(\psi_{x(x \wedge y)}\,c)$;
because the $\psi$ are bijections,
$\psi_{x(x \wedge y)}\,c = a$ and $c = \psi_{x(x \wedge
y)}^{-1}\,a$\trfootnote{The original has the typo
$\psi_{x(x \wedge y)}\,a$ for $\psi_{x(x \wedge y)}^{-1}\,a$.}.
Likewise, $c = \psi_{y(x \wedge y)}^{-1}\,b$.  Thus, (iii) follows from (ii).
The implication from (iv) to (i) is obtained similarly using
(20\supdelta).\trfootnote{The original has the typos ``The implication
from (iv) to (ii) is obtained similarly using (20).''}

We now consider the smallest equivalence relation on $L$, which equates
the range of each of the maps $\psi_{yx}$ with its domain.
To do this, we define a relation $\sim$ on $L$ as follows:
If $a \in K_x$
and $b \in K_y$, then $a \sim b$ should hold if and only if there is a
$z \in S$ with
$\psi_{zx}\,a = \psi_{zy}\,b$ or (equivalently by prop.~4.1)
$\psi_{xz}^{-1}\,a = \psi_{yz}^{-1}\,b$.
It is clear that this is exactly the equivalence relation we are
looking for:
\iitem{(21)}\textit{$\sim$ is an equivalence relation on $L$.}
\lastiitem

\textit{Proof.} Let $a_i \in K_{x_i}, (i = 0, 1, 2)$,
$a_0 \sim a_1$ and $a_1 \sim a_2$, i.e.\ by the proposition,
\begin{equation*}
\psi_{(x_0 \vee x_1)x_0}\,a_0 = \psi_{(x_0 \vee x_1)x_1}\,a_1
\textrm{\quad and\quad}
\psi_{(x_1 \vee x_2)x_1}\,a_1 = \psi_{(x_1 \vee x_2)x_2}\,a_2.
\end{equation*}

Since $x_1 \leqq x := (x_0 \vee x_1) \wedge (x_1 \vee x_2)$,
$\psi_{xx_1}\,a_1$, $\psi_{(x_0 \vee x_1)x}\,a_1$ and
$\psi_{(x_1 \vee x_2)x}\,a_1$ are by (19) defined, so by
(20\supdelta),\trfootnote{The original has the typo
(20) for (20\supdelta).}
also $\psi_{yx}\,a_1$ with $y := x_0 \vee x_1 \vee x_2$.
Now because of (19), $\psi_{y(x_0 \vee x_1)}
(\psi_{(x_0 \vee x_1)x_1}\,a_1)$ is defined and we have
$\psi_{yx_0}\,a_0 = \psi_{xx_1}\,a_1 = \psi_{yx_2}\,a_2$
and thus $a_0 \sim a_2$. Therefore $\sim$ is transitive.
Reflexivity and symmetry are trivially satisfied.

Now let $L^\pi := L/\sim$ be the quotient set and
$\pi: L \rightarrow L^\pi$ the canonical projection.%
\trfootnote{The original has the typo $K$ for $L$ in this sentence.}
For every $x \in S$ let the lattice $L_x^\pi$\trfootnote{The original
has the typo $L_x$ for $L_x^\pi$.}
be the homomorphic image
of $L_x$ under the map $\pi_x := (\pi|K_x) \circ i_x$.

\textbf{Theorem 4.2.} \textit{If the lattices $L_x\ (x \in S)$ and the
mappings $\phi_{yx}\ (x, y \in S,
x \leqq y)$ are an $S$-connected system, then the mappings
$\pi_x: L_x \rightarrow L^\pi_x$
for all $x \in S$ are isomorphisms and $L_x\ (x \in S)$ is an $S$-glued system.}

The $S$-glued sum of this system becomes the \textit{$S$-connected sum}
of the lattices $L_x$ under the identifications
$\phi_{yx}\ (x, y \in S, x \leqq y)$.
The $S$-glued system is monotone if and only if for all $x \prec y$,
$\im \phi_{yx} \neq L_y$.\bugfootnote{This condition must be updated due
to the updated definition of ``monotone'':
The $S$-glued system is monotone if and only if for all $x \prec y$,
$\im \phi_{yx} \neq L_y$ and $\dom \phi_{yx} \ne L_x$.}{2}

\textit{Proof.} If $a, b \in K_x$ and $\pi\,a = \pi\,b$, then there is a
$z \in S$ with $\psi_{zx}\,a = \psi_{zx}\,b$;
therefore $a=b$, since $\psi_{zx}$ is injective.
Thus $\pi | K_x$ is injective. There remain
conditions (1) to (4) from the definition of an $S$-glued system
to prove.

\hskip 0pt\\{[page 11/265]}\\\nopagebreak

Regarding (1) and (3): For $x \leqq y$ we have
$L^\pi_x \cap L^\pi_y = \pi \dom \psi_{yx} = \pi \im \psi_{yx}$ by
4.1.(ii), so by (17), (if it $\neq \emptyset$) it is a
filter in $L_x$ and an ideal in $L_y$; further
by (18), it is nonempty if $x \prec y$.

Regarding (2): Let $a, b \in L^\pi_x \cap L^\pi_y$ and
$a \leqq b$ in $L^\pi_x$. Then there are archetypes
$a^\prime, b^\prime$ of $a, b$ under $\pi$ in $K_x$ with $a^\prime \leqq b^\prime$
in $K_x$. On the other hand, there are archetypes
$a^{\prime\prime}, b^{\prime\prime}$ of $a, b$
under $\pi$ in $K_y$  with $a^{\prime\prime} \leqq b^{\prime\prime}$.
By 4.1.(ii),
$\psi_{(x \vee y)x}\,a^\prime = \psi_{(x \vee y)x}\,a^{\prime\prime}$
and $\psi_{(x \vee y)y}\,b^\prime = \psi_{(x \vee y)y}\,b^{\prime\prime}$.
Then $\psi_{(x \vee y)x}\,a^\prime \leqq \psi_{(x \vee y)x}\,b^\prime$
in $K_{x \vee y}$, and also
$a^{\prime\prime} = \psi^{-1}_{(x \vee y)y}(\psi_{(x \vee y)y}\,a^{\prime\prime})
\mathrel{\leqq\mskip-\thickmuskip\footnotemark}%
\footnotetext{The original omits $\leqq$.}
\psi^{-1}_{(x \vee y)y}(\psi_{(x \vee y)y}\,b^{\prime\prime}) = b^{\prime\prime}$
in $K_y$. Thus $a^{\prime\prime} \leqq b^{\prime\prime}$ in $K_y$
and $a \leqq b$ in $L^\pi_y$.

Regarding (4): If $a \in L^\pi_x \cap L^\pi_y$,
then there are $a^\prime \in K_x$ and $a^{\prime\prime} \in K_y$
with $\pi\,a^\prime = a = \pi\,a^{\prime\prime}$
and $\psi_{(x \vee y)x}\,a^\prime = \psi_{(x \vee y)y}\,a^{\prime\prime}$.
So $a = \pi(\psi_{(x \vee y)x}\,a^\prime) \in L^\pi_{x \vee
y}$. Using proposition~4.1,
$a \in L_{x \wedge y}$ can be proved similarly.

If $S$ is modular, there is a much simpler construction of $S$-glued sums:

Let S be a modular lattice of finite length, let $L_x\ (x \in S)$ be a family
of finite-length lattices and $\phi_{yx}\ (x, y \in S, x \prec y)$
a family of bijective maps.
The system consisting of the families $L_x\ (x \in S)$ and
$\phi_{yx}\ (x,y \in S, x \prec y)$ is called a \textit{locally $S$-connected
system} if for all
$u, v, x, y \in S$ with $u \prec v$ and
$x \wedge y \prec x, y \prec x \vee y$ the following apply:
\iitem{(22)}\textit{$\phi_{vu}$ is an isomorphism of a (nonempty) filter of
$L_u$ to an ideal of $L_v$.}
\iitem{(23)}$\phi_{(x \vee y)x} \circ \phi_{x(x \wedge y)} =
\phi_{(x \vee y)y} \circ \phi_{y(x \wedge y)}$
\iitem{(24)}$\im \phi_{(x \vee y)x} \cap \im \phi_{(x \vee y)y} \subseteq
\im \phi_{(x \vee y)(x \wedge y)}$
\iitem{(24\supdelta)}$\dom \phi_{x(x \wedge y)} \cap
\phi_{y(x \wedge y)} \subseteq \dom \phi_{(x \vee y)(x \wedge y)}$
\lastiitem

Here $\phi_{(x \vee y)(x \wedge y)}$ --- unique by (23) --- is defined as
\begin{equation*}
\phi_{(x \vee y)x} \circ \phi_{x(x \wedge y)}.
\end{equation*}

To get an $S$-connected system from this, we set
$\phi_{xx} = \id_{L_x}$ and choose for all $x < y$ a maximal chain
$x = x_0 \prec x_1 \prec x_2 \prec \cdots \prec x_{n-1} \prec y$.
We define
\begin{equation*}
\phi_{yx} := \phi_{yx_{n-1}} \circ \cdots \circ \phi_{x_2x_1} \circ \phi_{x_1x}.
\end{equation*}
This definition is independent of the choice of chain,
as can be shown by induction:
Let $x = y_0 \prec y_1 \prec \cdots \prec y_{m-1} \prec y$ be
another maximal chain. From modularity it follows that $n = m$. If $n = 1$,
there is nothing to show. So let $n > 1$. If $x_1 = y_1$, then by
the induction hypothesis, we have
$\phi_{yx_{n-1}} \circ \cdots \circ \phi_{x_2x_1} =
\phi_{yy_{n-1}} \circ \cdots \circ \phi_{y_2x_1}$
and the conclusion follows immediately.
So let $x_1 \neq y_1$. Set $z_2 = x_1 \vee y_1$; because of modularity
$z_2 \succ x_1, y_1$. Now choose a maximal chain
$z_2 \prec z_3 \prec \cdots \prec z_n = y$.
\\[\baselineskip]{[page 12/266]}\\[\baselineskip]\nopagebreak
Then by induction hypothesis it follows that
\begin{equation*}
\phi_{yx_{n-1}} \circ \cdots \circ \phi_{x_2x_1} =
\phi_{yz_{n-1}} \circ \cdots \circ \phi_{z_2x_1}\trfootnote{The
original has the typo $\phi_{z_2y_1}$ for $\phi_{z_2x_1}$.}
\end{equation*}
and
\begin{equation*}
\phi_{yy_{n-1}} \circ \cdots \circ \phi_{y_2y_1} =
\phi_{yz_{n-1}} \circ \cdots \circ \phi_{z_2y_1}
\end{equation*}
Because of (23), however, $\phi_{z_2x_1} \circ \phi_{x_1x} =
\phi_{z_2y_1} \circ \phi_{y_1x}$, and main statement follows immediately.

\textbf{Theorem 4.3.} \textit{It is assumed that $S$ is a modular
lattice of finite
length and the lattices $L_x\ (x \in S)$ and the mappings
$\phi_{yx}\ (x, y \in S, x \prec y)$
are a local $S$-connected system.
If the mappings $\phi_{yx}$ for all $x \leqq y$ in $S$ are defined as
above, the $L_x\ (x \in S)$ and $\phi_{yx}\ (x, y \in S, x \leqq y)$ form
an $S$-connected system, its $S$-connected sum is
called the locally $S$-connected sum of the $L_x\ (x \in S)$
and $\phi_{yx}\ (x, y \in S, x \prec y)$.
The associated $S$-connected system is monotone iff
$\im \phi_{yx} \neq L_y$ for all $x \prec y$.\bugfootnote{This
condition must be updated due
to the updated definition of ``monotone'':
The $S$-connected system is monotone if and only if for all $x \prec y$,
$\im \phi_{yx} \neq L_y$ and $\dom \phi_{yx} \ne L_x$.}{2}}

\textit{Proof.} The validity of (17), (18) and (19) is obvious.
Regarding (20). We use induction over the length $l$ of
$[x \wedge y, x \vee y]$.
If $l = 1$, then $x$ and $y$ are comparable and there is nothing to show.
If $l = 2$, then either $x$ and $y$ are comparable or (24) can be used
immediately. The final case is
$l \geqq 3$ and $\phi_{(x \vee y)x}\,a = \phi_{(x \vee y)y}\,b$ for
some $a \in L_x$, $b \in L_y$.
Then $x$ and $y$ cannot both
be lower neighbors of $x \vee y$, so there is a $z$
with $x < z < x$\trfootnote{The original has the typo $x \vee y$ for
$x$}.
Then $z \vee y = x \vee y$ and $x \wedge (z \vee y) = x \wedge y$;
also because of modularity,
$x \wedge y < z \wedge y$ and $x \vee (z \wedge y) = z$.
By (19), $\phi_{zx}\,a$ and $\phi_{(x \vee y)z} (\phi_{zx}\,a)$
exist.  By the induction hypothesis,
there is $c \in L_{z \wedge y}$ with $\phi_{(x \vee y)(z \wedge y)}\,c =
\phi_{(x \vee y)x}\,a$. Because of (19) and the injectivity of
$\phi_{(x \vee y)z}$
it follows that $\phi_{zx}\,a = \phi_{z(z \wedge y)}\,c$.
Applying the induction hypothesis again shows there is
$d \in L_{x \wedge y}$ with $\phi_{z(x \wedge y)}\,d
= \phi_{zx}\,a$ and thus $\phi_{(x \vee y)(x \wedge y)}\,d
= \phi_{(x \vee y)x}\,a$.
Dually, (20\supdelta) is proved using (24\supdelta).

\textit{Example.} Let $L_0, L_1$ be projective planes (that is, simple atomistic
modular lattices of length 3), $L_0$ generated by four lines $g_1, g_2,
g_3, g_4$, and $L_1$ generated by four points $P_1, P_2, P_3, P_4$.
Let $S$ and $L_{x_i}\ (i = 1, 2, 3, 4)$
be as in figure~1.\trfootnote{Figure~1 does not reproduce well.  The
second lattice is $L_{x_1}$ and all of its points (other than $e_0$)
have subscript $1$.  The third lattice is $L_j$ for $j = 2,3,4$ and
all of its points have subscript $j$.}
Define $\phi_{x_i0}$ by $\phi_{x_i0}\,g_i = 0_{x_i}$, $\phi_{x_i0}\,1_0 = a_i$; and
define $\phi_{1x_i}$ by $\phi_{1x_i}\,a_i =0_1$,
$\phi_{1x_i}\,1_{x_i} = P_i$; for $i = 1, 2, 3, 4$.
\\
\resizebox{7in}{!}{\includegraphics{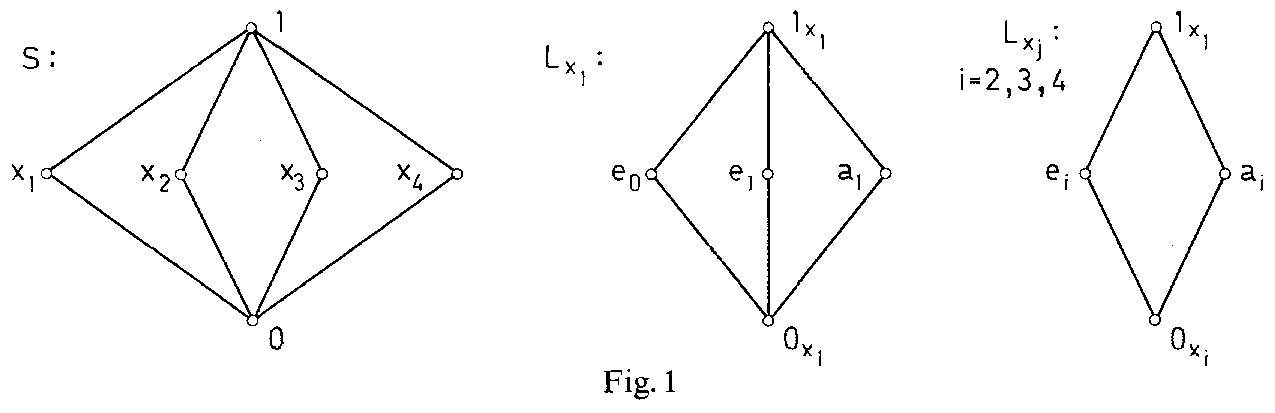}}

\hskip 0pt\\{[page 13/267]}\\\nopagebreak

This results in a locally $S$-connected system, the sum of which is a
simple modular lattice of length~6 and breadth~3 with generators
$e_0, e_1, e_2, e_3, e_4$.
The settles whether there exists a simple modular
lattice of any length with breadth $b = 3$ and 5~generators.
(The analogous problem with $b = 2$ and four generators was solved
positively by
Day, Wille and the author in [2].)

\vskip\baselineskip
\centerline{\textbf{5. $S$-gluing of homomorphisms}}
\nopagebreak\vskip\baselineskip

J\'onsson in [7] has developed the gluing of homomorphisms as a
counterpart to the gluing of lattices in Hall--Dilworth [4].
Again with restriction to finite-length lattices,
we make the following generalization:

\textbf{Theorem 5.1.} \textit{Let $L$ be the $S$-glued sum of lattices
$L_x\ (x \in S)$ and let
$\phi_x\ (x \in S)$ be a family of lattice homomorphisms
$\phi_x: L \rightarrow L^\prime$ such
that for all $x \prec y$ in $S$:
$\phi_x|L_x \cap L_y = \phi_y|L_x \cap L_y$.
Furthermore, for all $x, y$ in $S$, let:
\begin{equation}
\tag{*} \phi_x\,0_x + \phi_y\,0_y = \phi_{x \vee y}\,0_{x \vee y}
\textrm{\quad and\quad}
\phi_x\,1_x \cdot \phi_y\,1_y = \phi_{x \wedge y}\,1_{x \wedge y}.
\end{equation}
Then $\phi := \bigcup_{x \in S} \phi_x$ is a homomorphism from $L$
into $L^\prime$, which is injective if all $\phi_x\ (x \in S)$ are injective.}

\textit{Proof.}  We first show (without using (*)) that
$\phi_x|L_x \cap L_y = \phi_y|L_x \cap L_y$ for all $x, y \in S$.
We show this for $x \leqq y$ by induction on $y$.
The case $x = y$ is trivial and also
$y = x$\trfootnote{The original has the typo $y = 0$ for $y = x$.}
is the start of the induction.
Now let $x < y$. Then there is a $z \in S$ with
$x \leqq z \prec y$ and the
induction hypothesis gives $\phi_x|L_x \cap L_z = \phi_z|L_x \cap L_z$.
However, by (6),
$L_x \cap L_y \subseteq L_z$, so $\phi_x|L_x \cap L_y = \phi_z|L_x \cap L_y
= \phi_y|L_x \cap L_y$.
This shows the statement is true for all $x \leqq y$.
Now let $x, y \in S$ be arbitrary.
Then we already have shown
$\phi_x|L_x \cap L_{x \vee y} = \phi_{x \vee y}|L_x \cap L_{x \vee y}$
and $\phi_y|L_y \cap L_{x \vee y} = \phi_{x \vee y}|L_y \cap L_{x \vee y}$,
and by (4), $L_x \cap L_y \subseteq L_{x \wedge y}$.
A similar argument then shows
$\phi_x|L_x \cap L_y = \phi_{x \wedge y}|L_x \cap L_y = \phi_y|L_x \cap L_y$.

The proof that $\phi$ is a homomorphism is now carried out analogously
to the proof of theorem~2.1.
Let $a \in L_x$ and $y \geqq x$. We show $\phi(a + 0_y ) = \phi\,a
+ \phi\,0_y$ by induction on $y$.
If $x = y$ (which occurs in particular when $x = y = 0$), then
by (14), $a+0_x = a \uss+x 0_x$ therefore
$\phi(a+0_x) = \phi_x(a \uss+x 0_x) = \phi_x\,a + \phi_x\,0_x =
\phi\,a + \phi\,0_x$.
If $y>x$, then there is a $z$ with $x \leqq z \prec y$ and
\begin{equation*}
    a+0_y = a+0_z+0_y = (a+0_z) \uss+z 0_y.%
\trfootnote{The original has the typo $\uss+y$ for $\uss+z$.}
\end{equation*}
Therefore
\begin{equation*}
\phi(a + 0_y) = \phi_z((a + 0_ z) + 0_y) = \phi_z(a + 0_z) + \phi_z\,0_y
            = \phi\,a + \phi\,0_z + \phi_y\,0_y = \phi\,a+ \phi\,0_y.%
\trfootnote{The original has the typo $\phi_y$ for $\phi_z$.}
\end{equation*}
If $a \in L_x$ and $b \in L_y$ are arbitrary, then by (16),
\begin{equation*}
a + b = (a + 0_{x \vee y}) \uss+{x \vee y} (b + 0_{x \vee y}),
\end{equation*}
{[page 14/268]}\\[\baselineskip]\nopagebreak
therefore
\begin{equation*}
       \phi(a+b) = \phi_{x \vee y}(a+0_{x \vee y})+\phi_{x \vee y}(b+0_{x \vee y})
= \phi\,a + \phi\,0_{x \vee y} + \phi\,b
= \phi\,a + \phi\,0_x + \phi\,b + \phi\,0_y = \phi\,a + \phi\,b.
\end{equation*}
The dual follows similarly using ``$\cdot$''.

Let all $\phi_x$ be injective. 
For $a \prec b$ in $L$, by theorem~2.1 there is an $x \in S$ with
$a, b \in L_x$ and
$a \uss\prec x b$, from which $\phi_x\,a < \phi_x\,b$ and
$\phi\,a < \phi\,b$ result.
Since $L$ is finite length, this implies that if $a < b$ in $L$, then
$\phi\,a < \phi\,b$.
\emph{[Added by tr.:}
Now consider if $a$ and $b$ are incomparable in $L$, so that $a \vee b
> a, b$.  We have already
proven that $\phi$ is a homomorphism, so if $\phi\,a = \phi\,b$ then
$\phi(a \vee b) = \phi\,a = \phi\,b$, contradicting the preceding result.
Thus in all cases, if $a \neq b$, $\phi\,a \neq \phi\,b$.\emph{]}

\textbf{Addendum 5.2.} \textit{If $S$ is modular, then the requirement
\emph{(*)} is unnecessary.}

\textit{Proof.} The condition (*) has to be verified for all $x, y \in S$.
If $x \prec y$, then by (5), $0_y \in L_x$,
and $\phi_y\,0_y = \phi_x\,0_y \geqq \phi_x\,0_x$, implying the conclusion.
For $x \leqq y$ we argue similarly by induction.

If $x \wedge y \prec x, y \prec x \vee y$, then follows
$0_x, 0_y \in L_{x \wedge y}$%
\trfootnote{The original has the typo $L_{x \wedge y} \cap L_{x \vee y}$
for $L_{x \wedge y}$.}
by (3) and (5), so
\begin{equation*}
\phi_x\,0_x + \phi_y\,0_y = \phi_{x \wedge y}\,0_x \uss+{x \wedge y} \phi_{x \wedge y}\,0_x
= \phi_{x \wedge y}(0_x \uss+{x \wedge y}\,0_y)
= \phi_{x \wedge y}\,0_{x \vee y} = \phi_{x \vee y}\,0_{x \vee y}
\end{equation*}
by (14) and (13).

The general case is now proved by induction on the length $l$
of $[x \wedge y, x \vee y]$.
It remains to prove $l \geqq 3$:%
\trfootnote{The original has the typo $=$ for $\leqq$.}
let $z \in S$
with $x \wedge y < z < y$; by induction $\phi_x\,0_x + \phi_z\,0_z =
\phi_{x \vee z}\,0_{x \vee z}$ and
$\phi_{x \vee z}\,0_{x \vee z} + \phi_y\,0_y = \phi_{x \vee y}\,0_{x \vee y}$.
Therefore $\phi_{x \vee y}\,0_{x \vee y} =
\phi_x\,0_x + \phi_z\,0_z + \phi_y\,0_y$, and since $z \leqq y$, we
have already proven $\phi_z\,0_z \leqq \phi_y\,0_y$, so
$\phi_{x \vee y}\,0_{x \vee y} = \phi_x\,0_x + \phi_y\,0_y$.

The dual argument shows that
$\phi_x\,1_x \cdot \phi_y\,1_y = \phi_{x \wedge y}\,1_{x \wedge y}$.

\textbf{Corollary 5.3.} \textit{The $S$-glued sum of simple lattices is simple.}

\textbf{Corollary 5.4.} \textit{Let $L_x\ (x \in S$) be an $S$-glued system of
sublattices
of a lattice $L^\prime$. If $S$ is modular (or otherwise
$0_x + 0_y = 0_{x \vee y}$
and $1_x \cdot 1_y = 1_{x \wedge y}$ for all $x, y \in S$), then the
$S$-glued sum of
$L_x\ (x \in S)$ is a sublattice of $L^\prime$.}

\textbf{Corollary 5.5.} \textit{Let $L$ be the $S$-connected sum of
the system $L_x\ (x \in S)$,
$\phi_{yx}\ (x, y \in S, x \leqq y)$.  Let $\phi_x\ (x \in S)$ be a family
of homomorphisms $\phi_x: L_x \rightarrow L^\prime$
so that $\phi_x = \phi_y \circ \phi_{yx}$ for all $x \prec y$ from $S$.
Furthermore, let $S$ be modular or otherwise
for all $x, y \in S$, (*) is true.
For $\gamma_x := \phi_x \circ \pi_x^{-1}: L^\pi_x%
\trfootnote{The original has the typo $L_x$ for $L^\pi_x$.}
\rightarrow L^\prime\ (x \in S)$
then $\gamma = \bigcup_{x \in S} \gamma_x$ is a homomorphism of $L^\pi$ into
$L^\prime$, which is injective if and only if
all the $\phi_x$ are injective.}

\textit{Proof.}
\begin{align*}
\gamma_x|L^\pi_x \cap L^\pi_y
& = \phi_x \circ \pi_x^{-1}|\pi \dom \psi_{yx}
= \phi_x \circ \pi_x^{-1}|\pi\,i_x \dom \phi_{yx}
= \phi_x \circ \pi_x^{-1}|\pi_x \dom \phi_{yx}
= (\phi_x|\dom \phi_{yx}) \circ \pi_x^{-1}\footnotemark \\
& = \phi_y \circ \phi_{yx} \circ \pi_x^{-1} \\
& = (\phi_y|\im \phi_{yx}) \circ \pi_y^{-1}
= \phi_y \circ \pi_y^{-1}|\pi_y \im \phi_{yx}
= \phi_y \circ \pi_y^{-1}|\pi\,i_y \im \phi_{yx}
= \phi_y \circ \pi_y^{-1}|\pi \im \psi_{yx}
= \gamma_y|L_x^\pi \cap L_y^\pi
\end{align*}\footnotetext{The original omits ``$\circ\,\pi_x^{-1}$''
here.}
for all $x \prec y$, so the assumptions of theorem~5.1 are fulfilled.

\hskip 0pt\\{[page 15/269]}\\\nopagebreak

\vskip\baselineskip
\centerline{\textbf{6. The skeleton of a finite-length modular lattice}}
\nopagebreak\vskip\baselineskip

We recall the definition given in the introduction:
the \textit{skeleton} $S(M)$ of a finite-length modular lattice $M =
(M, +, \cdot)$ is the set of the smallest elements of the maximal
atomistic intervals of $M$.

Since it is known (cf. Wille [10; Thm.~1.1 and~1.4]) that a
finite-length modular lattice in which the largest element is a meet
of atoms is atomistic (and coatomistic), for each element a of M one
obtains atomistic intervals $[a,a^*]$ and $[a^+,a]$, where
\begin{equation*}
a^* = \begin{cases}
\sup \{b \mvert a \prec b\} & \text{if } a < 1 \\
1                           & \text{if } a = 1
\end{cases}
\end{equation*}
\begin{equation*}
a^+ = \begin{cases}
\inf \{b \mvert b \prec a\} & \text{if } a > 0 \\
0                           & \text{if } a = 0.
\end{cases}
\end{equation*}

\textbf{Lemma 6.1.} \textit{In a modular lattice $M$ of finite
length:\trfootnote{In cases (c) and (c\supdelta), the original has the
typo $\leqq$ for $=$.}
\begin{align*}
& \textrm{(a) If $a \leqq b$, then $a^* \leqq b^*$.} &
    & \textrm{(a\supdelta) If $a \leqq b$, then $a^+ \leqq b^+$.} \\
& \textrm{(b) $a \leqq a^{+*} \leqq a^*$.} &
    & \textrm{(b\supdelta) $a^+ \leqq a^{*+} \leqq a$.} \\
& \textrm{(c) $a^{*+*} = a^*$.} &
    & \textrm{(c\supdelta) $a^{+*+} = a^+$.} \\
& \textrm{(d) If $a = a^{*+}$ and $b = b^{*+}$, then $a + b = (a + b )^{*+}$.} &
    & \textrm{(d\supdelta) If $a = a^{+*}$ and $b = b^{+*}$,
        then $a \cdot b = (a \cdot b)^{+*}$.} \\
& \textrm{(e) $a^+ + b^+ = (a + b)^+$.} &
    & \textrm{(e\supdelta) $a^* \cdot b^* = (a \cdot b)^*$.}
\end{align*}}

\textit{Proof.} Regarding (a) and (a\supdelta) Since the conclusion is
trivial if $a = b$, let $a < b$.\trfootnote{The
original has the typo $=$ for $<$.} For all $c \succ a$,
modularity implies either $c \leqq b$ or $b \prec b+c$.
In either case, $c \leqq b^*$. Therefore $a^* \leqq b^*$.
Dually $a^+ \leqq b^+$.

Regarding (b): $a^+ \leqq a \leqq a^*$ is trivial, so by (a\supdelta),
$a^+ \leqq a^{*+}$.
Since the interval $[a,a^*]$ is coatomistic, the
intersection of its coatoms is $a$, so $a \geqq a^{*+}$.
So (b\supdelta) follows and, because of duality, also (b).

Regarding (c): By (b), $a^* \leqq a^{*+*}$; by (b\supdelta),
$a^{*+} \leqq a$ and then by (a),
$a^{*+*} \leqq a^*$; therefore $a^{*+*} = a^*$.
Dually, $a^{+*+} = a^+$.

Regarding (d): Since $a \leqq a + b$, it follows from (a) and
(a\supdelta) that
$a = a^{*+} \leqq (a + b)^{*+}$; likewise
$b \leqq (a + b)^{*+}$, so $a + b \leqq (a + b)^{*+}$.
By (b\supdelta), $(a + b)^{*+} \leqq a + b$, so
$(a + b)^{*+} = a + b$. The dual argument yields (d\supdelta).

We first prove equation (e) for the special case that
there is an element $c$ in $M$ with $c \prec a,b$.
Then $a^+, b^+ \leqq c$, so $a^+ \leqq a^+ + b^+ \leqq a$.
Therefore $[a^+ + b^+,a]$ is atomistic and $a \leqq (a^+ + b^+)^*$.
Likewise, $b \leqq (a^+ + b^+)^*$ and therefore
$a + b \leqq (a^+ + b^+)^*$. By (a\supdelta),
$(a+b)^+ \leqq (a^+ + b^+)^{*+}$.
By (c\supdelta) and (d), $(a^+ + b^+)^{*+} = a^{+*+} + b^{+*+}
= a^+ + b^+$. Hence $(a+b)^+ \leqq a^+ + b^+$.
On the other hand, it follows from (a\supdelta) that
$a^+, b^+ \leqq (a + b)^+$ and from this finally $(a + b)^+ = a^+ + b^+$.

The general case is proven by induction over $h(a) + h(b)$ ---
where $h(a)$ is the length of the interval $[0, a]$. If $h(a)+h(b)=0$, then
\\[\baselineskip]{[page 16/270]}\\[\baselineskip]\nopagebreak
$a=b=0$ and the statement becomes trivial. Let $h(a)+h(b)>0$ and without
loss of generality $a>0$. Then there is an $a_0 \prec a$. Applying
the induction hypothesis, $a_0^+ + b^+ = (a_0 + b)^+$.%
\trfootnote{The original has the typo $a_0^+ b^+$ for $a_0^+ + b^+$.}
If $a \leqq a_0 + b$, then
\begin{equation*}
a^+ + b^+ = a^+ + a_0^+ + b^+ = a^+ + (a_0+b)^+ = a^+ + (a+b)^+ = (a+b)^+.
\end{equation*}
So let $a \nleqq a_0 + b$.
Let $a_0 \prec a_1 \prec \cdots \prec a_m = a_0 + b$ be a maximal chain in $M$.
We show by induction on $i$: $a^+ + a_i = (a+a_i)^+$ for $i = 0, \ldots, m$.
This is trivial for $i=0$ and for $i>0$ follows from the induction hypothesis
and the already proven special case:
\begin{equation*}
a^+ + a_i^+ = a^+ + a_{i-1}^+ +a_i^+ = (a + a_{i-1})^+ + a_i^+ =
(a + a_{i-1} + a_i)^+ = (a + a_i)^+,
\end{equation*}
because $a_{i-1} \prec a_i$ and $a_{i-1} \prec a_{i-1} + a$. Further,
\begin{equation*}
a^+ + b^+ = a^+ + a_0^+ + b^+ = a^+ + (a_0+b)^+ = a^+ +a_m^+ =
(a+a_m)^+ = (a+b)^+.\trfootnote{The original has the typo $a_n$ for $a_m$.}
\end{equation*}
This ends the proof of (e) and --- because of duality --- also of (e\supdelta).

\textbf{Theorem 6.2.}
\textit{Let M be a modular lattice of finite length with largest
element 1 and smallest element 0. Then for the set $S(M)$ of the
smallest elements of the maximal atomistic intervals of $M$:}

\vskip0.5\baselineskip
(f)\quad \textit{$S(M) = \{x \mvert x \in M \textrm{ and }
x^{*+} = x \}$ and the maximal atomistic
intervals of $M$ are exactly the intervals $[x, x^*]$ with $x \in S(M)$
(equivalently, with $x^{*+} = x$).}

\vskip0.5\baselineskip
(g)\quad \textit{$S(M)$ is a sub-join-semilattice of $M$ with
largest element
$1^+$ and smallest element 0. The lattice structure of the skeleton
$S(M)=(S(M), \vee, \wedge)$ is given by
$x \vee y = x + y$ and $x \wedge y = (x \cdot y)^{*+}$.}
\vskip0.5\baselineskip

\textit{Proof.} Regarding (f): Let $[a,b]$ be a maximal atomistic interval.
$[a, b]$ is contained in the atomistic interval $[a, a^*]$,
so $b = a^*$. Dually
it follows that $a = b^+$ and therefore $a=a^{*+}$.
Conversely, let $a = a^{*+}$. The interval
$[a, a^*]$ is atomistic, so there is a maximal atomistic interval
$[c, d]$ with $[c, d] \supseteq [a, a^*]$. Then $[c, a^*]$ is also atomistic, so
$a = a^{*+} \leqq c$ and therefore $a = c$,
which implies $d \leqq a^*$, so $d = a^*$ and $a \in S(M)$.

Regarding (g): $x \vee y = x + y$ follows immediately from (d), $0 = 0^{*+}$
from (b\supdelta); further because of
(a\supdelta), $a^{*+} \leqq 1^+$ for all $a \in M$.
From (c), $(x \cdot y)^{*+} \in S(M)$.
From (b\supdelta) we get $(x \cdot y)^{*+} \leqq x,y$.
For all $z \in S(M)$ for which $z \leqq x, y$, because of (a) and (a\supdelta),
$z = z^{*+} \leqq (x \cdot y)^{*+}$.
So $(x \cdot y)^{*+} = x \wedge y$.

The preliminary work of sections~1, 2, 5 and~6 now makes possible:

\textit{Proof of the main theorem.} By theorem~6.2, the
maximal atomistic
intervals of $M$ are exactly the intervals $[x, x^*]$ with $x \in S(M)$.
It must be shown that the lattices $[x, x^*]\ (x \in S(M))$ are
an $S(M)$-glued
system and fulfill the requirements of corollary~5.4. Then
the $S$-glued sum of the $[x, x^*]\ (x \in S(M))$ is a sublattice
$\bigcup_{x \in S(M)} [x, x^*]$ of $(M, +, \cdot)$ and,
because $M$ is finite-length, $M = \bigcup_{x \in S(M)} [x, x^*]$
is identical to $(M, +, \cdot)$.

\hskip 0pt\\{[page 17/271]}\\\nopagebreak

By theorem~6.2, $S(M)$ is a lattice and, since it is contained in $M$
as a partial order, it is of finite length;
furthermore $x \vee y = x + y$ and
$(x \wedge y)^* = (x \cdot y)^{*+*} = (x \cdot y)^* = x^* \cdot y^*$
by (e\supdelta).
(1) and (2) are trivially fulfilled.
Regarding (4): Let $a \in [x, x^*] \cap [y,y^*]$.
Then $x \vee y = x + y \leqq a \leqq x^* \leqq (x \vee y)^*$
and $x \wedge y = (x \cdot y)^{*+} \leqq x \cdot y \leqq a \leqq
x^* \cdot y^* = (x \wedge y)^*$. Finally, (3) follows
from the following

\textbf{Lemma 6.3.} \textit{If $x \prec y$ in $S(M)$,
then $x < y \leqq x^*$ in $M$.}

\textit{Proof.} Choose $a$ such that $x \leqq a \prec y$ in $M$.
Then $y \leqq a^*$.
Necessarily, $a^{*+} \in S(M)$ and $x = x^{*+} \leqq a^{*+} \leqq y^{*+} = y$,
and since $x \prec y$ in $S(M)$, either $a^{*+} = x$ or $a^{*+} = y$.
But $a^{*+} \leqq a \prec y$, so $a^{*+} = x$ and $a^* = x^*$.
Therefore $y \leqq x^*$.

Dual to the skeleton is the \textit{dual skeleton $S^\delta(M)$}, consisting of
the largest elements of the maximal atomistic intervals.
Then $S^\delta(M) = \{ a \mvert a \in M \textrm{ and } a^{+*} = a \}$
is a sub-meet-semilattice of $M$.

\textbf{Theorem 6.4.} \textit{The skeleton and the dual skeleton of a
finite-length modular lattice
are isomorphic via the mutually inverse isomorphisms $*$ and $+$.}

The proof is a trivial consequence of (a) and (c) of lemma~6.1.

\textbf{Theorem 6.5.} \textit{The dual skeleton of a finite-length
modular lattice is the dual of%
\footnote{The original omits ``dual of''.}
the skeleton of the dual lattice.}

This is immediately clear, since every atomistic interval is coatomistic
and vice versa.

\textbf{Corollary 6.6.} \textit{The skeleton of a self-dual
finite-length modular lattice is self-dual.}

Of the applications of the main theorem, we mention only only this:
the breadth or $n$-distributivity of a modular lattice can be interpreted
as an ``inherent dimension''.
This is expressed in the the following two consequences.
(Here ``projective geometry'' stands for
``lattice of subspaces of a projective geometry''):

\textbf{Consequence 6.7.} \textit{For a finite-length modular lattice the
following are equivalent:
\begin{enumerate}
\item[\emph{(i)}] $M$ has breadth $\leqq n$.
\item[\emph{(ii)}] No (possibly reducible) projective geometry
of dimension $n$ can be embedded in $M$.
\item[\emph{(iii)}] All $[x, x^*]\ (x \in S(M))$ are (possibly reducible) projective
geometries of dimension $\leqq n-1$.
\end{enumerate}}

\textbf{Consequence 6.8.} \textit{For a finite-length modular lattice
the following are equivalent:
\begin{enumerate}
\item[\emph{(i)}] $M$ is $n$-distributive.
\end{enumerate}}
\vskip\baselineskip{[page 18/272]}\\\nopagebreak   
\textit{\begin{enumerate}
\item[\emph{(ii)}] No irreducible projective geometry of dimension $n$
can be embedded in $M$.
\item[\emph{(iii)}] For every $x \in S(M)$, $[x, x^*]$ is a direct product of
irreducible projective geometries of dimension $\leqq n-1$.
\end{enumerate}}

\textit{Proof.} Consequence~6.7 is a direct consequence of theorem~3.3
and the theorem of Birkhoff [1; Chap.~IV, \S~7] referenced in the
introduction.

Regarding 6.8: (i) $\Rightarrow$ (ii) and (i) $\Rightarrow$ (iii)
follow directly from the result given in section~3
from Huhn [5; theorem~1.1]. (ii) $\Rightarrow$ (iii) is
trivial and (iii) $\Rightarrow$ (i) follows from theorem~3.4.
(i) $\Rightarrow$ (ii) was proved by Huhn more generally for
for algebraic modular lattices (letter communication).

\vskip\baselineskip
\centerline{\textbf{7. The class of all skeletons}}
\nopagebreak\vskip\baselineskip

The skeleton of a finite length modular lattice is in general not
modular. On the contrary, theorem~7.2 shows that even the class of
skeletons of distributive lattices generates the class of
all lattices. (By J\'onsson [8; Cor.~5.2], the class of all
lattices is generated by the finite lattices).

\textbf{Lemma 7.1.} \textit{If $M$ is the $S$-glued sum of modular
lattices $M_x\ (x \in S)$, then $S(M)= \bigcup_{x \in S} S(M_x)$.
$S$ is isomorphic to $S(M)$ by the map
$x \mapsto 0_x$ and $M_x = [0_x, 0_x^*]$ for all $x \in S$
if and only if $M_x\ (x \in S)$ is a monotone $S$-glued system of
atomistic modular lattices.}%
\bugfootnote{The lemma is not correct without further
conditions, but assuming that the system is the stronger form of
``monotone'', the other clauses of the lemma follow.
The uses of lemma~7.1 in this article all involve strongly monotone
systems, and so are unaffected by this problem.}{3}

\textit{Proof.} By theorem~2.1, for all $x \in S$ the intervals of
$M_x$ are also intervals of $M$.
Because of lemma~3.1 there is, for every atomistic
interval $[a, b]$ of $M$, an $x \in S$ such that $[a, b] \subseteq M_x$.
Therefore an interval is a
maximally atomistic interval in $M$ if and only if this is the case in an $M_x$
for a suitable $x \in S$. Thus $S(M)= \bigcup_{x \in S} S(M_x)$.
The rest of the claim
now follows directly from theorem~2.1 and lemma~6.3.

\textbf{Theorem 7.2.} \textit{Every finite lattice is isomorphic to
the skeleton of a finite distributive lattice.}

\textit{Proof.} Let $S$ be a finite lattice. Let $\power(S)$ be the power set
lattice of $S$ and $\power(S)^\delta$ be its dual lattice.
For every $x \in S$ let
$(x]$ be the ideal generated by $x$ in S and $[x)$ be the filter
generated by $x$ in $S$.
Let $L_x$ be the interval $[(\emptyset, [x)), ((x], \emptyset)]$ in
$\power(S) \times \power(S)^\delta\ (x \in S)$.%
\footnote{Perhaps more intuitively stated as
$L_x = \power((x]) \times \power([x))^\delta$.}
The $L_x\ (x \in S)$ form an $S$-glued system.
For $x \prec y$, $(\{y\}, \{y\}) \in L_y \setminus L_x$,\bugfootnote{To
prove the stricter definition of ``monotone'', we must also note that
$(\{x\}, \{x\}) \in L_x \setminus L_y$.}{2} and thus
the $S$-glued system is monotone.
Since the $L_x$ are themselves distributive atomistic lattices,
it follows from theorem~3.4 and lemma~7.1 that the $S$-glued sum of
$L_x\ (x \in S)$ is distributive and its skeleton is isomorphic to $S$.

\textbf{Theorem 7.3.} \textit{Every finite-length modular lattice $S$
is isomorphic to the skeleton of a sublattice of $S \times S$.}

\hskip 0pt\\{[page 19/273]}\\\nopagebreak

\textit{Proof.} For every $x \in S$, let $L_x$ be the interval
$[(x^+, x), (x, x^*)]$ of $S \times S$.
By lemmas 6.1 and 6.3 it is easy to see that $L_x\ (x \in S)$ is a monotone
$S$-glued system of atomistic lattices.
If $M$ is the $S$-glued sum,
then $S$ is isomorphic to $S(M)$ by lemma~7.1. Furthermore
by 6.1(e)
\begin{equation*}
0_x + 0_y = (x^+, x) + (y^+, y) = ((x+y)^+, x+y)
= ((x \vee y)^+, x \vee y) = 0_{x \vee y}
\end{equation*}
and similarly $1_x \cdot 1_y = 1_{x \wedge y}$. So $M$ is a sublattice of
$S \times S$ by corollary 5.4.

Every finite chain is a skeleton --- e.g.\ of the chain with one
more element. But since every non-trivial ultraproduct of finite
chains has infinite length, the class of all skeletons%
\footnote{Note that ``skeletons'' implicitly means ``skeletons of
finite-length lattices'', and so all skeletons are themselves
finite-length.}
is certainly not axiomatizable. However,

\textbf{Theorem 7.4.}  \textit{Let $n$ be a natural number.
Then the class of all skeletons of modular lattices of length $n$
(length $\leqq n$) is an axiomatic class.%
\trfootnote{Here, ``axiomatic class'' means a class axiomatizable by
first-order sentences, not just equations.
Also, ``elementary sublattice'' means a sublattice of a
lattice that is an elementary substructure of that lattice.  An
elementary substructure of a structure is one that satisfies exactly
the same first-order sentences as the structure does.
See [T3; Thm.~2.13].}}

\textit{Proof.}
It is sufficient to show that the class is closed under taking of
of ultraproducts and elementary sublattices.
Regarding taking ultraproducts:
Let $M_i\ (i \in I)$ be modular lattices of length $n$ (lengths $\leqq n$),
$\scrF$ be an ultrafilter on $I$, $S = \uss\bigtimes{i \in I^\scrF} S (M_i)$ and
$M = \uss\bigtimes{i \in I^\scrF} M_i$.
$M$ is a modular lattice of length $n$ (length $\leqq n$).
Define $\phi$ to map each element of $S$
--- an equivalence class of elements of $\uss\bigtimes{i \in I} S(M_i)$ ---
containing $(x_i \mvert i \in I)$ (with $x_i \in S(M_i) \subseteq M_i$)
to the element of $M$
--- an equivalence class of elements of $\uss\bigtimes{i \in I} M_i$ ---
containing that same $(x_i \mvert i \in I)$.
By the structure of the ultraproduct, $\phi$ is a well-defined map from $S$
to $M$ and is an isomorphism.
Due to theorem~6.2, the property $x \in S(M)$ is equivalent to
$x = x^{*+}$, and so is first-order.
This implies that since each $S(M_i)$ is the skeleton of $M_i$,
$\phi(S)$ is the skeleton of $M$.

Regarding taking elementary sublattices:
Let $M$ a modular lattice of length $n$ (length $\leqq n$) and
$S^\prime$ an elementary sublattice of $S = S(M)$.
Then $x \prec y$ in $S^\prime$ implies $x \prec y$ in $S$.
Therefore the atomistic intervals $[x,x^*]\ (x \in S^\prime)$ form a
monotone $S^\prime$-glued system whose sum, $U = \bigcup_{x \in S^\prime} [x, x^*]$,
is a sublattice of $M$ (by corollary~5.4) and
$S^\prime$ is isomorphic to --- even equal to --- $S(U)$ (by lemma~7.1).
Furthermore, $S^\prime$ and $S$ have a maximum chain in common.
From the main theorem and theorem~2.1 it follows that any refinement of this
chain is already in $U$; so $U$ also has length $n$ (length $\leqq n$).
Thus $S^\prime$ is the skeleton of $U$, a lattice of length $n$
(length $\leqq n$).

\vskip\baselineskip
\centerline{\textbf{References}%
\trfootnote{The items with the format ``[Xxxx9999x]'' are bibliography
entry codes in \url{https://alum.mit.edu/www/worley/Math}.}}
\nopagebreak\vskip\baselineskip

\begin{enumerate}
\item[1.] Birkhoff, G.: Lattice Theory, third edition. Providence:
Amer.\ Math.\ Soc.\ 1967.  [Birk1967a]
\item[2.] Day, A, Herrmann, C., Wille, R.: On modular lattices with four generators. Preprint 1972, submitted to Algebra Universalis
\item[3.] Gr\"atzer, G., Lakser, H., J\'onsson, B.: The Amalgamation Property in equational classes of modular lattices. Preprint 1971. Summary in Notices Amer.\ Math.\ Soc.\ \textbf{18}, \#635 (1971).  [GratzJ\'onsLaks1973a]

\hskip 0pt\\{[page 20/274]}\\\nopagebreak

\item[4.] Hall, M., Dilworth, R. P.: The imbedding problem for modular lattices. Ann.\ of Math.\ \textbf{45}, 450--456 (1944).  [HallDil1944a]
\item[5.] Huhn, A.: Schwach distributive Verb\"ande.  [Weakly distributive lattices.]  Acta F.R.N.\ Univ.\ Comen. --- Mathematica --- mimoriadne \v{c}islo --- 51--56 (1971).
\item[6.] J\'onsson, B.: Modular lattices and Desargues' Theorem. Math.\ Scandinav.\ \textbf{2}, 293--314 (1954).  [J\'ons1954a]
\item[7.] J\'onsson, B.: Arguesian lattices of dimension $n \leqq 4$. Math.\ Scandinav.\ 7, 133 - 145 (1959).  [J\'ons1959b]
\item[8.] J\'onsson, B.: Algebras whose congruence lattices are distributive. Math.\ Scandinav.\ \textbf{21}, 110--121 (1967).  [J\'ons1967a]
\item[9.] Khee-Meng Koh: On the Frattini sublattice of a lattice. Algebra Universalis \textbf{1}, 104--116 (1971).
\item[10.] Wille, R.: Halbkomplement\"are Verb\"ande. [Semi-complementary lattices.]   Math.\ Z.\ \textbf{94}, 1--31 (1966).
\item[T1.] Herrmann, C.: $S$-verklebte Summen von Verb\"anden.  Math.\ Z. \textbf{130}, 255--274 (1973).  [Herr1973a]
\item[T2.] Day, A., Herrmann, C.: Gluings of modular lattices. Order \textbf{5}, 85--101 (1988).  [DayHerr1988a]
\item[T3.] Frayne, T., Morel, A., Scott, D: Reduced direct products. Fund.\ Math.\ \textbf{51}, 195--228 (1962).{\\}[FrayneMorScottD1962a]
\item[T4.] Bandelt, H.-J.: Tolerance relations on lattices.  Bull.\ Austral.\ Math.\ Soc. \textbf{23}, 367--381 (1981).  [Band1981a]
\item[T5.] Day, A.\ and Freese, R.: The Role of
Gluing Constructions in Modular Lattice Theory, in The Dilworth
theorems: selected papers of Robert P. Dilworth (Kenneth P. Bogard, Ralph
S. Freese, and Joseph P. S. Kung, eds.), Contemporary Mathematicians,
Springer, New York,
1990, 251--260.  [DayFrees1990a]
\item[T6.] Dilworth, R.\ P., The arithmetical theory of Birkhoﬀ
lattices, in The Dilworth theorems: selected papers of
Robert P. Dilworth (Kenneth P. Bogard, Ralph S.
Freese, and Joseph P. S. Kung, eds.), Contemporary Mathematicians,
Springer, New York, 1990, 286--289.  [Dil1941a]

\end{enumerate}
\vskip\baselineskip

\begin{samepage}
\noindent
Dr.\ C.\ Herrmann\\
Mathematisches Institut\\
der Technischen Hochschule\\
D-6100 Darmstadt\\
Bundesrepublik Deutschland\\
\\
\\
(Received September 28, 1972)

\end{samepage}

\newpage

\begin{center}
\textbf{TRANSLATOR'S NOTES} \\[\baselineskip]
{\small DALE R. WORLEY}
\end{center}
\vskip 2\baselineskip

\textbf{1. $S$-glued sums are finite length.}[Satz~2.1, Th.~2.1]\enspace
$S$-glued sums need not be finite-length lattices.
A counterexample to Theorem~2.1 can be constructed as an elaboration
of the straightforward gluing shown in figure~\ref{fig:3by3}.  The
lattice constructed in figures~\ref{fig:infinite-length-parts}
and~\ref{fig:infinite-length-sum} has chains of all finite lengths.
To prove finite length, additional conditions are needed.
If the skeleton is finite length and the sum is modular, that forces
its elements to be rankable,
which suffices to show it is finite-length.  If the component lattices
are modular, the sum is modular.

\begin{figure}[b]
\hbox to \linewidth{\hfill
\begin{tikzcd}[cramped,sep=small,every arrow/.append style={dash}]
          & 4 \ar[dl]\ar[dr] \\
2 \ar[dr] &                  & 3 \ar[dl] \\
          & 1 \\
	  & \textup{Skeleton}
\end{tikzcd}
\hfill
\begin{tikzcd}[cramped,sep=small,every arrow/.append style={dash}]
          & e \ar[dl]\ar[dr] \\
b \ar[dr] &                  & c \ar[dl] \\
          & a \\
          & L_1
\end{tikzcd}
\hfill
\begin{tikzcd}[cramped,sep=small,every arrow/.append style={dash}]
          & g \ar[dl]\ar[dr] \\
d \ar[dr] &                  & e \ar[dl] \\
          & b \\
          & L_2
\end{tikzcd}
\hfill
\begin{tikzcd}[cramped,sep=small,every arrow/.append style={dash}]
          & h \ar[dl]\ar[dr] \\
e \ar[dr] &                  & f \ar[dl] \\
          & c \\
          & L_3
\end{tikzcd}
\hfill
\begin{tikzcd}[cramped,sep=small,every arrow/.append style={dash}]
          & i \ar[dl]\ar[dr] \\
g \ar[dr] &                  & h \ar[dl] \\
          & e \\
          & L_4
\end{tikzcd}
\hfill
\begin{tikzcd}[cramped,sep=small,every arrow/.append style={dash}]
          &                      & i \ar[dl]\ar[dr] \\
          & g \ar[dl]\ar[dr]     &                  & h \ar[dl]\ar[dr] \\
d \ar[dr] &                      & e \ar[dl]\ar[dr] &                  & f \ar[dl] \\
          & b \ar[dr]            &                  & c \ar[dl] \\
          &                      & a \\
	  &                      & \textup{Sum}
\end{tikzcd}
\hfill}
\caption{A sum lattice of finite length: skeleton, components, and sum}
\label{fig:3by3}
\end{figure}

\begin{figuregroup}

\begin{figure}[b]
\hbox to \linewidth{\hfill
\begin{tikzcd}[cramped,sep=small,every arrow/.append style={dash}]
            & \Omega \ar[dl]\ar[d]\ar[dr] \\
1 \ar[dr] & 2 \ar[d]                  & \cdots \ar[dl] \\
            & A \\
	    & \textup{Skeleton}
\end{tikzcd}
\hfill
\begin{tikzcd}[cramped,sep=small,every arrow/.append style={dash}]
            & e \ar[dl]\ar[d]\ar[dr] \\
b_1 \ar[dr] & b_2 \ar[d]                  & \cdots \ar[dl] \\
            & a \\
	    & L_A
\end{tikzcd}
\hfill
\begin{tikzcd}[cramped,sep=small,every arrow/.append style={dash}]
            & g_n \ar[dl]\ar[ddr] \\
d_{nn} \ar[d]\\
\vdots \ar[d] &                    & e \ar[dddl] \\
d_{n2} \ar[d] \\
d_{n1} \ar[dr] \\
            & b_n \\
	    & L_n
\end{tikzcd}
\hfill
\begin{tikzcd}[cramped,sep=small,every arrow/.append style={dash}]
            & h \ar[dl]\ar[d]\ar[dr] \\
g_1 \ar[dr] & g_2 \ar[d]                  & \cdots \ar[dl] \\
            & e \\
	    & L_\Omega
\end{tikzcd}
\hfill}
\caption{A sum lattice with chains of unbounded length: skeleton and components}
\label{fig:infinite-length-parts}
\end{figure}

\begin{figure}[b]
\hbox to \linewidth{\hfill
\begin{tikzcd}[cramped,sep=small,every arrow/.append style={dash}]
              &                     &                                          & h\ar[dll]\ar[d]\ar[drr]\ar[drrrr] \\
              & g_1\ar[ddl]\ar[ddr] &                                          & g_2\ar[ddl]\ar[dr] &               & g_3\ar[ddlll,bend right=25]\ar[dr] &             & \cdots \\
              &                     &                                          &                    & d_{22}\ar[dd] &                                    & d_{33}\ar[d] \\
d_{11}\ar[ddr] &                    & e\ar[ddl]\ar[ddr]\ar[ddrrr,bend right=25] &                   &               &                                    & d_{32}\ar[d] & \cdots \\
              &                     &                                          &                   & d_{21}\ar[dl] &                                     & d_{31}\ar[dl] \\
              & b_1\ar[drr]         &                                          & b_2\ar[d]          &              & b_3\ar[dll]                         &              & \cdots\ar[dllll] \\
              &                     &                                          & a
\end{tikzcd}
\hfill}
\caption{A sum lattice with chains of unbounded length: sum}
\label{fig:infinite-length-sum}
\end{figure}

\end{figuregroup}

\textbf{2. When 0 elements are unequal.}[(12)]\label{p:0-unequal}\enspace
Contrary to statement (12), $0_x$ can be equal to $0_y$ even if $x < y$ and
$L_y \not\subset L_x$.  In particular, consider $x = 1$ and $y = 2$ in the
sum lattice of figure~\ref{fig:overlap}.

\begin{figure}[b]
\hbox to \linewidth{\hfill
\begin{tikzcd}[cramped,sep=small,every arrow/.append style={dash}]
4 \ar[d] \\
3 \ar[d] \\
2 \ar[d] \\
1 \\
\textup{Skeleton}
\end{tikzcd}
\hfill
\begin{tikzcd}[cramped,sep=small,every arrow/.append style={dash}]
b \ar[d] \\
a \\
L_1
\end{tikzcd}
\hfill
\begin{tikzcd}[cramped,sep=small,every arrow/.append style={dash}]
d \ar[d] \\
c \ar[d] \\
b \ar[d] \\
a \\
L_2
\end{tikzcd}
\hfill
\begin{tikzcd}[cramped,sep=small,every arrow/.append style={dash}]
d \ar[d] \\
c \\
L_3
\end{tikzcd}
\hfill
\begin{tikzcd}[cramped,sep=small,every arrow/.append style={dash}]
d \ar[d] \\
c \\
L_4
\end{tikzcd}
\hfill
\begin{tikzcd}[cramped,sep=small,every arrow/.append style={dash}]
d \ar[d] \\
c \ar[d] \\
b \ar[d] \\
a \\
\textup{Sum}
\end{tikzcd}
\hfill}
\caption{Overlapping of components: skeleton, components, and sum}
\label{fig:overlap}
\end{figure}

The fix for this problem is a stricter definition of ``monotone'':

\textbf{Definition T.1}
\textit{An $S$-glued lattice system is \emph{monotone} (in the stricter sense)
if for all
$x \prec y$, $x, y \in S$, then $L_y \not\subset L_x$ and
$L_x \not\subset L_y$.}

Note this renders the definition of ``monotone'' self-dual.
The first clause of this definition (given in the
original) ensures the map $x \mapsto 1_x$ is an injection;
the second clause of this definition ensures the
map $x \mapsto 0_x$ is an injection.

\textbf{3. When a gluing is not monotone.}[Lem.~7.1]\enspace
Lemma~7.1 is not correct without further conditions.  A
counterexample is figure~\ref{fig:monotone}; the skeleton of the sum
is $\{0\}$, but $\bigcup_{x \in S} S(M_x) = S(L_1) \cup S(L_2) = \{0, b\}$.
However, if the $S$-gluing is monotone in the stricter, self-dual
sense, then all of the clauses of the lemma follow.
All uses of the lemma in the article involve monotone gluings,
so the remainder of the article is unaffected by this problem.

\begin{figure}[b]
\hbox to \linewidth{\hfill
\begin{tikzcd}[cramped,sep=small,every arrow/.append style={dash}]
2 \ar[d] \\
1 \\
\textup{Skeleton}
\end{tikzcd}
\hfill
\begin{tikzcd}[cramped,sep=small,every arrow/.append style={dash}]
                 & 1 \ar[dl]\ar[d]\ar[dr] \\
ab \ar[d]\ar[dr] & ac \ar[dl]\ar[dr]      & bc \ar[dl]\ar[d] \\
a \ar[dr]        & b \ar[d]               & c \ar[dl] \\
	         & 0 \\
                 & L_1
\end{tikzcd}
\hfill
\begin{tikzcd}[cramped,sep=small,every arrow/.append style={dash}]
                 & 1 \ar[dl]\ar[dr] \\
ab \ar[dr]       &                        & bc \ar[dl] \\
                 & b \\
                 & L_2
\end{tikzcd}
\hfill
\begin{tikzcd}[cramped,sep=small,every arrow/.append style={dash}]
                 & 1 \ar[dl]\ar[d]\ar[dr] \\
ab \ar[d]\ar[dr] & ac \ar[dl]\ar[dr]      & bc \ar[dl]\ar[d] \\
a \ar[dr]        & b \ar[d]               & c \ar[dl] \\
	         & 0 \\
                 & \textup{Sum}
\end{tikzcd}
\hfill}
\caption{When a gluing is not monotone: skeleton, components, and sum}
\label{fig:monotone}
\end{figure}

Lemma~7.1 can be split into these parts:

\textbf{Lemma T.2}
\textit{If $M$ is the $S$-glued sum of modular
lattices $M_x\ (x \in S)$,
then $S(M) \subset \bigcup_{x \in S} S(M_x)$
and $S^\delta(M) \subset \bigcup_{x \in S} S^\delta(M_x)$}
\begin{proof}
\leavevmode

By theorem~2.1, for all $x \in S$ the intervals of
$M_x$ are also intervals of $M$.
Because of Lemma~3.1 there is, for every atomistic
interval $[a, b]$ of $M$, an $x \in S$ such that $[a, b] \subset M_x$.
Therefore a
maximally atomistic interval in $M$ is a maximally atomistic interval
in an $M_x$
for a suitable $x \in S$. Thus $S(M) \subset \bigcup_{x \in S} S(M_x)$.
Dually,  $S^\delta(M) \subset \bigcup_{x \in S} S^\delta(M_x)$.
\end{proof}

\textbf{Lemma T.3}
\textit{If $M$ is the strongly monotone $S$-glued sum of atomistic modular lattices
$M_x\ (x \in S)$, then the $M_x$ are exactly the set of maximal
atomistic intervals of M.}
\begin{proof}
\leavevmode

For every $x \in S$, $M_x$ is an atomistic interval in $M$.
If $M_x$ was contained in a larger atomistic
interval $L$, by Lemma~3.1, there is some $M_y$ that contains $L$, and
necessarily $y \neq x$.
Since $M_x \subset L \subset M_y$, by (7),
$M_x \subset M_x \cap M_y = M_{x \vee y}$.
Then $M_x$ is contained in $M_{x \vee y}$ and $x \vee y \neq x$, which
is impossible in a strongly monotone S-gluing.
Thus every $M_x$ is a maximal atomistic interval in $M$.

Conversely, let $L$ be a maximal atomistic interval in $M$.
By Lemma~3.1, there is some $x \in S$ such that $M_x$ contains $L$.
Since $L$ is maximal in $M$, $M_x = L$.
Thus, $L = M_x$ for some $x \in S$.
\end{proof}

\textbf{Lemma T.4}
\textit{If $M$ is a finite-length modular lattice, then $M$ is the $S$-glued
sum of the strongly monotone $S$-glued system which is the maximal atomistic
intervals of $M$, $[x, x^*]\ (x \in S(M))$.}
\begin{proof}
\leavevmode

That $M$ is the $S$-glued sum of the maximal atomistic intervals of
$M$ is the Main Theorem.  What remains to be proved is that the sum is
strongly monotone.
If $x \prec y$ in $S(M)$, then Lemma~6.3 shows that $x < y \leq x^*$,
which implies that $x \in [x, x^*] \setminus [y, y^*]$ and so
$[x, x^*] \not\subset [y, y^*]$.
To show that $y^* \in [y, y^*] \setminus [x, x^*]$ and so
$[y, y^*] \not\subset [x, x^*]$:  Assume that $y^* \in [x, x^*]$.
Then $y^* \leq x^*$ and by Lemma~6.1, $y^{*+} \leq x^{*+}$.
By Theorem~6.2, $x^{*+} = x$ and $y^{*+} = y$, so $y \leq x$,
which contradicts $x \prec y$.
\end{proof}

The remaining properties of $M$ given in Lemma~7.1 follow directly
from these lemmas.

\end{document}